%% file: main.tex
\documentclass{article}
\usepackage[a4paper, total={6.5in, 10in}]{geometry}
\input{preamble}
\usepackage{lscape}
\usepackage{authblk}
\usepackage{tabularray}

\title{Bridging the gap between agent based models and continuous opinion dynamics}
\author[1,2]{Andrew Nugent}
\author[2]{Susana N. Gomes}
\author[2]{Marie-Therese Wolfram}
\affil[1]{MathSys CDT, University of Warwick}
\affil[2]{Mathematics Institute, University of Warwick}

\begin{document}

\maketitle

\begin{abstract}
There is a rich literature on microscopic models for opinion dynamics; most of them fall into one of two categories - agent-based models or differential equation models - with a general understanding that the two are connected in certain scaling limits. In this paper we show rigorously this is indeed the case.
In particular we show that DEMs can be obtained from ABMs by simultaneously rescaling time and the distance an agent updates their opinion after an interaction. This approach provides a pathway to analyse much more diverse modelling paradigms, for example: the motivation behind several possible multiplicative noise terms in stochastic differential equation models; the connection between selection noise and the mollification of the discontinuous bounded confidence interaction function; and how the method for selecting interacting pairs can determine the normalisation in the corresponding differential equation. Our computational experiments confirm our findings, showing excellent agreement of solutions to the two classes of models in a variety of settings.
\end{abstract}


\textbf{Highlights:}
\renewcommand{\labelitemi}{{\tiny$\bullet$}}
\begin{itemize}
\itemsep-0.5em
        \item Proves the connection between agent-based and differential equation models of opinion formation. 
        \item Explains the relationship between selection noise and mollifying interaction functions.
        \item Motivates multiplicative diffusion in SDE models using additional noise terms in an ABM. 
        \item Demonstrates the convergence of each model through numerical simulations. 
\end{itemize}


\section{Introduction} \label{Section: Introduction}

Much of opinion dynamics finds its roots in the models of Hegselmann-Krause (HK) \cite{hegselmann2004opinion} and Deffuant-Weisbuch (DW) \cite{deffuant2000mixing}. Both models are based on the underlying assumption, known as bounded confidence, that people interact only with those who already share a sufficiently similar opinion. In the DW model each individual can only interact with one randomly chosen partner at a time, while in the HK model individuals can interact with everyone in their interaction radius simultaneously. Computational experiments have shown that solutions to both models exhibit similar macroscopic behaviours, such as a transition from consensus to a growing number of opinion clusters as the bounded confidence radius is reduced, but their exact dynamics can be quite different. 

Since their initial introduction, both models have been studied and adapted extensively, with their influence and major features visible in a wide range of modern opinion dynamics models (see review papers \cite{lorenz2007continuous,flache2017models,sirbu2017opinion,noorazar2020classical}). 

In this paper we will investigate the two most common approaches in opinion dynamics: agent-based models (ABMs) and differential equation models (DEMs). The category of ABMs includes the original DW model, while many DEMs are generalisations of the ordinary differential equation (ODE) version of the HK model. In ABMs opinion updates are discrete events in which a pair of randomly chosen individuals interact and then change their opinion(s). By contrast, in DEMs opinions change continually as individuals interact constantly with the entire population, rather than with one individual at a time. This distinction is more than the difference between discrete and continuous time, but is a fundamental difference in the nature of interactions and opinion updates. For example, in ABMs both the order of these random interactions and the distance by which individuals update their opinions can play a major role in the dynamics, but these features are essentially absent from DEMs. Despite these differences, we will show that the models are connected in certain scaling regimes. 

Our analysis begins with the following observation: reducing the update distance in the DW model leads to a slower convergence to equilibrium \cite{deffuant2000mixing,urbig2007communication} and affects both the location and composition of opinion clusters \cite{urbig2007communication,laguna2004minorities}. We take this situation to the extreme, considering an ABM in which agents interact at a very fast rate but with a very small update distance. In this setting, the effect of many small, random interactions in the ABM mirrors the continuous opinion updates characteristic of DEMs. Motivated by this, we will show that under a simultaneous rescaling of time and update distance, the trajectories of the ABM converge to those of a limiting DEM. 

Primarily this clarifies when ABMs and DEMs will give similar results, but also addresses the question of when the two models are capturing the same real-world behaviour, highlighting their inherent assumptions. Furthermore, we will show how changes made in the ABM can translate to those in the corresponding DEM, drawing parallels between the choices made in each setting. This link also bridges the gap between two research communities by demonstrating the essential similarity of two seemingly disparate modelling approaches. Our effort to establish connections between the many opinion formation models will help motivate modelling choices and improve our understanding of these complex processes, we hope this drives others to seek similar connections.

This paper is organised as follows: firstly in Section \ref{Section: Theory} we define an ABM and state our main results about its convergence to an ODE system. We also discuss several adaptations to the ABM that appear in the relevant literature, including the introduction of various additional noise terms, and give the corresponding ODE and stochastic differential equation (SDE) limits. A summary of the limiting models for all cases examined can be found in Table \ref{tab:Summary table} and proofs given in the Supplementary Material \ref{Supplementary Material}. In Section \ref{Section: Numerical Results} we present numerical results for each of the ODE and SDE limits described in Section \ref{Section: Theory}. Finally we conclude with a discussion of the significance and limitations of our results in Section \ref{Section: Discussion}.

\section{Theory} \label{Section: Theory}

We begin by introducing several ABMs, based on a generalisation of the DW model, and discuss their respective limits. Interested readers can find the rigorous derivation and the mathematical statements of these limits, which are based on Durrett \cite{durrett2018stochastic}, in the Supplementary Material. 

Throughout this paper we will always consider a finite population size $N$, and do not examine any mean-field limits or density-based models (see for example \cite{lorenz2007continuous,fennell2021generalized,chu2022density,ben2003bifurcations,fortunato2005vector,goddard2022noisy,como2009scaling}). Such models offer one way of comparing the behaviour of ABMs and DEMs (through the limit of a large population size) while the aim of this paper is to offer an alternative approach for a fixed $N < \infty$. 

\subsection{Agent based model} \label{Subsection: ABM}

Consider $N$ agents with initial opinions $x(0) = (x_1,\dots,x_N) \in \mathds{R}^N$. For generality, we do not prescribe an interval in which an individual's opinion should lie, but for the majority of models we consider opinions will remain within the convex hull of the initial opinions \cite{motsch2014heterophilious,ceragioli2021generalized}.

Beginning at $t=0$, we proceed in discrete timesteps of size $h>0$. At each timestep, choose a pair of individuals $i$ and $j$ from the population. This choice is made uniformly at random with replacement (different selection mechanisms are considered in Section \ref{Subsection: Normalisation}). After individuals $i$ and $j$ have been selected, they go on to interact with probability $p_{ij}(x)$. Note that $p_{ij}$ does not include the probability of selecting individuals $i$ and $j$, only the probability of interaction after their selection.

The form of $p_{ij}(x)$ plays a significant role in the consensus dynamics. The majority of opinion formation models assume that interactions depend on the opinion distance, that is $p_{ij}(x) = \phi(|x_i - x_j|)$ for some interaction function $\phi(d):\mathds{R}^+ \rightarrow [0,1]$. The bounded confidence (BC) interaction function with radius $R\in\mathds{R}^+$, as used in the DW model, also falls into this category,
\begin{equation} \label{Eqn: BC interaction function}
    \phi_R(d) = 
    \begin{cases}
        1 & \text{ if } d \leq R \,,\\
        0 & \text{ otherwise} \,.
    \end{cases}
\end{equation}
Here the only source of noise is the selection of individuals $i$ and $j$, as an interaction will occur (or not occur) deterministically based on their opinion difference. Other possible choices of interaction functions include exponentially decaying functions and polynomials, see e.g. \cite{como2009scaling,koponen2022agent,motsch2014heterophilious}. 

\begin{figure}[ht!]
    \centering
    \includegraphics[width=0.8\linewidth]{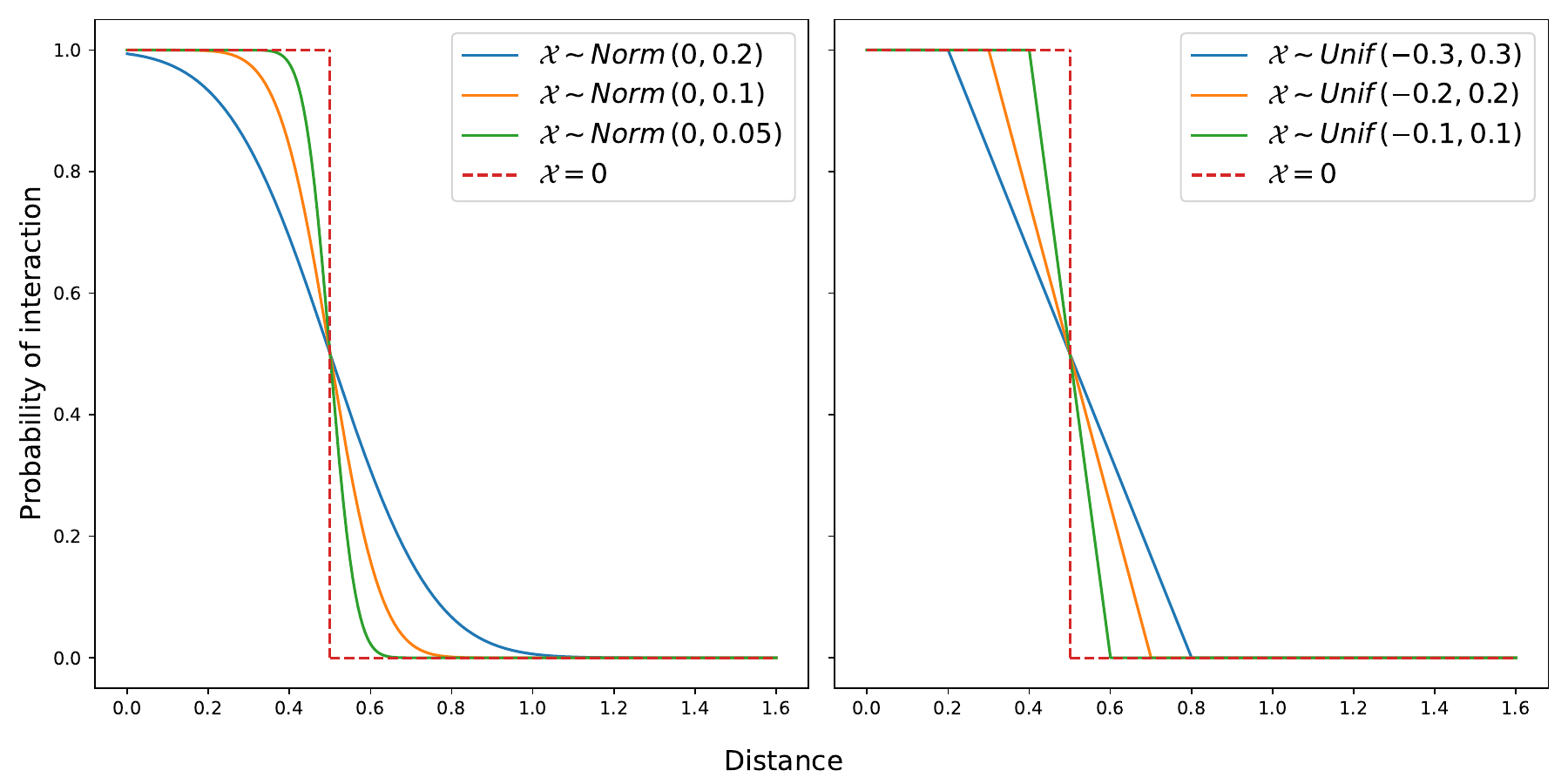}
    \caption{Examples of interaction functions of the form \eqref{Eqn: nosiy interaction function}, in which the discontinuous bounded confidence interaction function \eqref{Eqn: BC interaction function} (shown in red, dashed) is mollified with various types of selection noise. All the resulting interaction functions are Lipschitz continuous, but not necessarily differentiable everywhere.}
    \label{fig: example 'noisy' interaction functions}
\end{figure}

An common adaptation to the BC interaction function is the introduction of selection noise (see \cite{steiglechner2023noise} and references therein). Hereby, after picking $i$ and $j$, a random number $r$ is drawn from a given distribution $\mathcal{X}$. Individuals $i$ and $j$ then interact if
\begin{equation}
    |x_i - x_j| < R + r.
\end{equation}
Hence the probability of interaction is 
\begin{align}
    p_{ij}(x) 
    &= 1 - F_\mathcal{X}(|x_i - x_j| - R) \,, \label{Eqn: nosiy interaction function}
\end{align}
where $F_\mathcal{X}$ is the cumulative distribution function of $\mathcal{X}$. Thus selection noise corresponds to an altered interaction function. In particular, the introduction of selection noise can provide a mollification (or smoothing) of the discontinuous BC interaction function. Several examples are shown in Figure \ref{fig: example 'noisy' interaction functions} for a variety of normal and uniformly distributed $\mathcal{X}$'s. Note that for the normal distributions the interaction function is always strictly positive, while for the uniform distributions it is zero outside some interval around $R$. This distinction affects which interactions are possible and has a major impact on the behaviour of the limiting ODE systm \cite{motsch2014heterophilious}. 

The probability of interaction $p_{ij}$ also provides an opportunity to introduce a social network by setting the probability of interaction equal to the network adjacency matrix.

In order to establish the connection between ABMs and DEMs, we will require the following assumptions on $p_{ij}$:
\begin{assumption} \label{Assumption group: p_ij}
For all $i,j\in\{1,\dots,N\}$, $p_{ij}$ satisfies 
\begin{enumerate}[label=\alph*)]
    \item$p_{ij}(x)\in[0,1]$ for all $x\in\mathds{R}^N$. \label{Assumption: p_ij in 0,1}
    \item $p_{ij}(x)$ is globally Lipschitz continuous in $x$.  \label{Assumption: p_ij cts}
\end{enumerate}
\end{assumption}

Assumption \ref{Assumption group: p_ij}\ref{Assumption: p_ij in 0,1} is required so that $p_{ij}$ can be interpreted as the probability of $i$ and $j$ interacting once selected. Assumption \ref{Assumption group: p_ij}\ref{Assumption: p_ij cts} is required so that the limiting system of differential equations is well-defined on $\mathds{R}^N$. In particular this excludes the discontinuous BC interaction function without mollification by selection noise. Such assumptions on the continuity of interaction functions are often necessary to make progress in analysing DEMs \cite{ceragioli2012continuous,brooks2022emergence}, making the identification of realistic mechanisms, such as selection noise, that give rise to smoothed interaction functions especially important. 

If $i$ and $j$ do interact, then individual $i$ updates their opinion according to,
\begin{align} \label{eqn: opinion update rule}
    x_i(t+h) &= x_i(t) + \mu^h\, \big(x_j(t) - x_i(t) \big),
\end{align}
with all other individuals' opinions remaining the same in the next timestep. For a given timestep $h$, let $\mu^h = Nh$. If $h\in(0,N^{-1}]$ then $\mu \leq 1$ and individuals' opinions will remain within the convex hull of the initial opinions $x(0)$.

\begin{remark}
    The setup in which both individuals $i$ and $j$ update their opinions will give the same limit with $\mu = Nh/2$, under the assumption that $p_{ij}(x)=p_{ji}(x)$ for all $x\in\mathds{R}^N$.
\end{remark}

Overall, an interaction between individuals $i$ and $j$ gives the following possible updates to the opinion of individual $i$,
\begin{equation} \label{Eqn: Standard update scheme}
    x_i(t+h) = 
    \begin{cases}
        x_i(t) + \mu^h\, \big(x_j(t) - x_i(t) \big) & \text{ with probability } p_{ij}(x) \\
        x_i(t) & \text{ with probability } 1 - p_{ij}(x) \,.
    \end{cases}
\end{equation}

This process of selecting a pair of individuals, determining if they interact and updating the first individual's opinion, continues until some given time $T\in[0,\infty)$ is reached. 

\textbf{Main result:} Let $x_0$ be the opinions at time $t=0$. Then, as $h\rightarrow0$, trajectories of the ABM obtained using update rule \eqref{Eqn: Standard update scheme} converge in distribution to the solution $X_i = (X_i(t)\,:t\in[0,T])$ satisfying,
\begin{equation} \label{Eqn: ODE model}
    \frac{dX_i}{dt} = \frac{1}{N} \sum_{j=1}^N p_{ij}(X)\, (X_j - X_i) \,,\quad X(0) = x_0.
\end{equation}
Convergence in distribution means that the distribution of a sequence of random variables (here the trajectories of the ABM) converges to the distribution of a limiting random variable. Here the limit is the solution of the ODE system \eqref{Eqn: ODE model}, so it is deterministic, and thus the trajectories of the ABM also converge in probability. Informally, this means that as $h\rightarrow0$ the trajectories of the ABM are almost surely arbitrarily close to the solution of the ODE system. Formal definitions of both types of convergence are given along with the proof of the main result in the Supplementary Material.

The ODE \eqref{Eqn: ODE model} is a general model of opinion formation in continuous-time, whose behaviour is analysed in e.g. \cite{motsch2014heterophilious,blondel2010continuous,lacker2018mean}. This result provides a link between the ABM and ODE model, showing that under the assumption of small, frequent pairwise interactions their behaviours are comparable. In addition, it shows how the interaction probability $p_{ij}$ in the ABM, which can encode interaction functions, selection noise and networks, corresponds to opinion weights in the ODE model. We will continue to explore this connection by considering various relevant generalisations in the following sections. 

\begin{remark}
    Here we have defined an ABM that is continuous in opinion space and discrete in time. It would also be possible to define a continuous-time version of this ABM, with pairwise interactions occurring randomly at some constant rate multiplied by $p_{ij} h$. In this continuous-time ABM, interactions (and therefore opinion updates) are still discrete `events', unlike the continual opinion updates in DEMs. Hence, changing the ABM from discrete- to continuous-time would not resolve the fundamental difference between ABMs and DEMs. However, the re-scaling and convergence results used in this paper could also be applied to a continuous-time version of the ABM to obtain the same limiting models. 
\end{remark}

\subsection{Normalisation} \label{Subsection: Normalisation}

A common variation of the ODE model \eqref{Eqn: ODE model} is to replace the normalisation by the population size $N$ with normalisation by either the network degree (e.g. in \cite{nugent2023evolving}) or by the sum of interactions (e.g. in \cite{boghosian2022particle,motsch2014heterophilious,koponen2022agent}). We show below how both of these variations can be obtained as the limit of an ABM by changing the way in which the interacting individuals $i$ and $j$ are selected. Recall that in \eqref{Eqn: Standard update scheme}, individuals $i$ and $j$ are chosen uniformly at random with replacement. 

\begin{remark} \label{Remark: Choice w/o replacement, updating both opinions}
    If the choice is made without replacement, the ODE \eqref{Eqn: ODE model} is unchanged if we set $\mu = (N-1)h$, rather than $\mu = Nh$. 
\end{remark}

Firstly we consider normalisation by node degree. Define a fixed, weighted network with adjacency matrix $A\in[0,1]^{N \times N}$, with $A_{ii} = 1$ for all $i=1,\dots,N$. This network will only be used when picking the pair of individuals $i$ and $j$, and we do not assume that this network appears in the interaction probabilities $p_{ij}$. For each individual $i$, denote their node degree by $k_i$, this is 
\begin{equation}
    k_i = \sum\limits_{\ell=1}^N A_{i \ell} \,.
\end{equation}
At each timestep, select an individual $i$. The second individual $j$ is then selected with probability 
\begin{equation} \label{Eqn: node degree prob}
    \pi_j = \frac{A_{ij}}{k_i}.
\end{equation}
That is, $j$ is selected from amongst the neighbours of individual $i$, with probability proportional to the weight of the edge between $i$ and $j$. The limiting ODE model in this case is given by 
\begin{equation} \label{Eqn: ODE model, k_i normalisation}
    \frac{dX_i}{dt} = \frac{1}{k_i} \sum_{j=1}^N A_{ij}\,p_{ij}(X)\, (X_j - X_i).
\end{equation}
Secondly we consider normalisation by interaction probability. In order for this setup to be well-defined we assume that there exists $c>0$ such that $p_{ii}(X) > c$ for all $X\in\mathds{R}^N$. This can be interpreted as each individual always maintaining some confidence in their own opinion. 

At each timestep, first select an individual $i$. For this individual, calculate $p_{ij}$ for each $j\in\{1,\dots,N\}$. The second individual $j$ is then selected with probability 
\begin{equation} \label{Eqn: interaction prob normalisation}
    \pi_j(X) = \frac{p_{ij}(X)}{\sum\limits_{\ell=1}^N p_{i \ell}(X) }.
\end{equation}
That is, $j$ is selected with probability proportional to the likelihood that they would go on to interact with individual $i$. In this case we do not consider again the interaction probability $p_{ij}$, instead after individual $j$ is selected an interaction always then occurs. The limiting ODE model in this case is given by 
\begin{equation} \label{Eqn: ODE model, p_ij normalisation}
    \frac{dX_i}{dt} = \Bigg(\sum\limits_{\ell=1}^N p_{i \ell}(X)\Bigg)^{-1} \sum_{j=1}^N p_{ij}(X)\, (X_j - X_i).
\end{equation}
This more closely reflects the normalisation by the number of interacting agents in the original HK model. The key difference is that the ODE system \eqref{Eqn: ODE model, p_ij normalisation} corresponds to an ABM in which individual $i$ considers the likelihood of interaction with all individuals before $j$ is chosen, whereas the ODE system \eqref{Eqn: ODE model} corresponds to an ABM in which individual $i$ considers the likelihood of interaction only after $j$ is chosen. These variations show how the interpretation of the method through which $i$ and $j$ are selected in the ABM can help motivate the choice of normalisation in the ODE. In the case of opinion dynamics, normalisation by the entire population size may in fact be the more realistic option. 

\subsection{Additional sources of noise} \label{Section: Additional sources of noise}

As described in \cite{steiglechner2023noise} there are various places that additional noise can be introduced into the ABM. For many of these additional noise terms, as well as some not considered in \cite{steiglechner2023noise}, a similar re-scaling as in the previous section again gives rise to a DEM. The following cases will be considered:
\renewcommand{\labelitemi}{{$\bullet$}}
\begin{itemize}
    \item Ambiguity noise: in which the opinion of individual $j$ is not clearly communicated to individual $i$, meaning that $x_j$ is replaced by $\omega_j = x_j + \eta^h$, where $\eta^h$ is some random variable depending on $h$. 
    \item External noise: in which the opinion update of individual $i$ is altered by adding a random variable $\xi^h$. Alternatively, this noise could be added only when individuals $i$ and $j$ interact, in which case it is referred to as adaptation noise \cite{steiglechner2023noise}. 
    \item Random update distance: in which the fixed value of the update distance $\mu$ is replaced by a random variable $\nu^h$. 
\end{itemize}
These changes are made by altering the update rule \eqref{Eqn: Standard update scheme}, rather than by changing the way in which $i$ and $j$ are selected. We return to the original setup in which $i$ and $j$ are chosen uniformly at random with replacement. In this section we will introduce assumptions on the distributions of $\eta^h$, $\xi^h$ and $\nu^h$, and their behaviours as $h\rightarrow0$, to ensure convergence of the model with the same drift term as previously. We show that, under these assumptions, ambiguity noise leads to the same ODE limit \eqref{Eqn: ODE model}, while external noise and random update distance have the potential to give rise to SDEs. 

\begin{remark}
    Due to the many possible implementations, we do not consider here any boundary conditions on either the ABM or SDE system. Such boundary conditions could be included in the ABM and the changes followed through to determine the corresponding conditions on the SDE, but such an adaptation is beyond the scope of this paper. This is discussed further in Section \ref{Section: Discussion}.
\end{remark}

The following notation will be used throughout this section. Consider a family of real-valued random variables $\zeta = (\zeta^h)_{h>0}$, where $\zeta^h$ indicates dependence on $h$. We define
\begin{equation} \label{Eqn: m_k}
    m_k(\zeta) = \lim\limits_{h\rightarrow0} \dfrac{\mathds{E}\big[ (\zeta^h)^k \big]}{h},
\end{equation}
with $m_k(\zeta)$ undefined if a finite limit does not exist. 

\subsubsection{Ambiguity noise} \label{Subsection: ambiguity noise}

Ambiguity noise represents the fact that an individual may not clearly communicate their opinion. The noise therefore affects both the probability of interaction and the location of individual $i$'s new opinion. 

We assume here that $p_{ij}(x) = \phi(|x_j - x_i|)$ for some function $\phi:\mathds{R}^+\rightarrow[0,1]$ that is Lipschitz continuous with Lipschitz constant $L$. We modify the ABM by replacing the update rule \eqref{Eqn: Standard update scheme} with
\begin{equation} \label{Eqn: Ambiguity noise update scheme}
    x_i(t+h) = 
    \begin{cases}
        x_i(t) + \mu^h\, \big(\omega_j(t) - x_i(t) \big) & \text{ with probability } \phi(|\omega_j - x_i|) \\
        x_i(t) & \text{ with probability } 1 - \phi(|\omega_j - x_i|) \,,
    \end{cases}
\end{equation}
where $\omega_j = x_j + \eta^h$, with the following assumptions on $\eta$.
\begin{assumption} \label{Assumption group: eta^h}
Let $\eta^h$ correspond to the ambiguity noise in \eqref{Eqn: Ambiguity noise update scheme}. Assume the family of real random variables $\eta = (\eta^h)_{h>0}$ satisfies
\begin{enumerate}[label=\alph*)]
    \item $\lim\limits_{h\rightarrow0} \mathds{E}\big[ |\eta^h| \big] = 0$.  \label{Assumption: eta^h abs mean}
    \item There exists a constant $C\in(0,\infty)$, independent of $h$, such that $\mathds{E}\big[ (\eta^h)^2 \big] < C$ for all $h>0$. \label{Assumption: eta^h finite variance}
\end{enumerate}
\end{assumption} 
Broadly speaking, the first assumption ensures that $\phi(|\omega_j - x_i|)\rightarrow\phi(|x_j - x_i|)$ as $h\rightarrow0$, while the second is required to ensure convergence of the ABM. Under these assumptions on $\eta$, this system has the same ODE limit \eqref{Eqn: ODE model} as the ABM without ambiguity noise. Ambiguity noise affects the likelihood of interactions, but the first part of Assumption \ref{Assumption group: eta^h} ensures that in the limit $h\rightarrow0$ there is no net increase or decrease in this likelihood. As $\mu$ becomes small, each individual interaction has less significance, with the average behaviour instead driving the dynamics. Even though ambiguity noise does affect the ABM \cite{steiglechner2023noise}, this additional stochasticity does not impact the average behaviour and therefore does not lead to a different limit. This demonstrates that introducing additional noise and complexity at the ABM level does not necessarily impact the re-scaled system.

\subsubsection{External noise} \label{Subsection: External noise}

We now adapt the ABM by introducing additive noise to represent random external influences on an individual's opinion. We replace the update rule \eqref{Eqn: Standard update scheme} with
\begin{equation} \label{Eqn: External noise update scheme}
    x_i(t+h) = 
    \begin{cases}
        x_i(t) + \mu^h\, \big(x_j(t) - x_i(t) \big) + \xi^h & \text{ with probability } p_{ij}(x) \\
        x_i(t) + \xi^h & \text{ with probability } 1 - p_{ij}(x) \,,
    \end{cases}
\end{equation}
where $\xi^h\in\mathds{R}$ is a random variable representing all external influences on individual $i$'s opinion. Motivated by the results applied in the Supplementary Material, we introduce the following assumptions on $\xi$,
\begin{assumption} \label{Assumption group: xi^h}
Let $\xi^h$ correspond to the ambiguity noise in \eqref{Eqn: External noise update scheme}. Assume the family of real random variables $\xi = (\xi^h)_{h>0}$ satisfies
\begin{enumerate}[label=\alph*)]
    \item $m_1(\xi)\, = 0$.  \label{Assumption: xi^h mean}
    \item $m_2(\xi)$ exists. \label{Assumption: xi^h variance limit}
    \item If $m_2(\xi) > 0$, then $m_3(\xi) = m_4(\xi) = 0$. \label{Assumption: xi^h higher moments}
\end{enumerate}
\end{assumption}
The first of these assumptions ensures that the noise $\xi$ does not introduce a drift in opinion. 

If $m_2(\xi) > 0$, meaning the variance of $\xi$ is decreasing $o(h)$, we obtain an SDE limit with additive noise. An example of this is $\xi^h \sim \mathcal{N}(0,sh)$ for some constant $s>0$, which gives $m_2(\xi) = s$. The final assumption on $\xi$ is required to ensure convergence in this case. The limiting SDE is given by
\begin{equation} \label{Eqn: limiting SDE external noise}
    dX_i = \frac{1}{N} \sum_{j=1}^N p_{ij}(X)\, (X_j - X_i) \,\text{d}t + \bigg(\frac{m_2(\xi)}{N}\bigg)^{\frac{1}{2}} \,\text{d}\beta_i \,,
\end{equation}
where $\beta_i$ are independent standard Brownian motions. This shows that additive noise introduced in the ABM can translate directly to additive noise in the scaling limit. The noise is additive as $\frac{m_2(\xi)}{N}$ is a constant. Note that for $m_2(\xi) = 0$ we recover the ODE limit \eqref{Eqn: ODE model}.

\subsubsection{Adaptation noise}

Next we consider an adaptation of the ABM in which noise is added only if individuals $i$ and $j$ interact. This is referred to as adaptation noise in \cite{steiglechner2023noise}. We replace the update scheme in \eqref{Eqn: Standard update scheme} with
\begin{equation} \label{Eqn: External noise with interactions update scheme}
    x_i(t+h) = 
    \begin{cases}
        x_i(t) + \mu^h\, \big(x_j(t) - x_i(t) \big) + \xi^h & \text{ with probability } p_{ij}(x) \\
        x_i(t) & \text{ with probability } 1 - p_{ij}(x) \,,
    \end{cases}
\end{equation}
where $\xi$ is again assumed to satisfy Assumption \ref{Assumption group: xi^h}. In this case we obtain the limiting SDE,
\begin{equation} \label{Eqn: limiting SDE external noise on interactions}
    dX_i = \frac{1}{N} \sum_{j=1}^N p_{ij}(X)\, (X_j - X_i) \,\text{d}t + \Bigg(\frac{m_2(\xi)}{N^2} \sum_{j=1}^N p_{ij}(X)\Bigg)^{\frac{1}{2}}\,\text{d}\beta_i \,.
\end{equation}
Here the additional noise term in the ABM has the potential to change the structure of the limiting DEM, motivating a feasible multiplicative diffusion function for the SDE.

\subsubsection{Random update distance} \label{Subsection: Random update distance}

Finally we adapt the ABM by adding a new source of noise: replacing the fixed value of $\mu$ with a random variable $\nu$. This represents the idea that opinion updates do not move by a set difference, but that the outcome of an interaction is also random. At a particular interaction, individual $j$ may be more or less persuasive and individual $i$ may be more or less receptive or confident in their existing opinion. Additionally, we do not exclude the possibility of $\nu^h$ taking negative values, only requiring that the mean returns the fixed value of $\mu = N h$ used in the standard ABM, meaning individuals may move apart after an interaction. We replace the update scheme \eqref{Eqn: Standard update scheme} with
\begin{equation} \label{Eqn: Noisy update distance update scheme}
    x_i(t+h) = 
    \begin{cases}
        x_i(t) + \nu^h\, \big(x_j(t) - x_i(t) \big) & \text{ with probability } p_{ij}(x) \\
        x_i(t) & \text{ with probability } 1 - p_{ij}(x) \,,
    \end{cases}
\end{equation}
Unlike in the external noise case, the noise term is also multiplied by $(x_j - x_i)$, so a consensus state in which all individuals have identical opinions is still a fixed point.

We make the following assumptions on $\nu$, mirroring those on $\eta$,
\begin{assumption} \label{Assumption group: nu^h}
Let $\nu^h$ correspond to the ambiguity noise in \eqref{Eqn: Noisy update distance update scheme}. Assume the family of real random variables $\nu = (\nu^h)_{h>0}$ satisfies
\begin{enumerate}[label=\alph*)]
    \item $m_1(\nu)\, = N$.  \label{Assumption: nu^h mean}
    \item $m_2(\nu)$ exists. \label{Assumption: nu^h variance limit}
    \item If $m_2(\nu) > 0$, then $m_4(\nu) = 0$. \label{Assumption: nu^h higher moments}
\end{enumerate}
\end{assumption}

As before, if $m_2(\nu) = 0$ we obtain the standard deterministic ODE limit \eqref{Eqn: ODE model}. If $m_2(\mu)>0$ we obtain a limiting SDE with a new multiplicative diffusion function
\begin{equation} \label{Eqn: limiting SDE noisy update distance}
    dX_i = \frac{1}{N} \sum_{j=1}^N p_{ij}(X)\, (X_j - X_i) \,\text{d}t + \Bigg( \frac{m_2(\nu)}{N^2} \sum_{j=1}^N p_{ij}(X)\, (X_j - X_i)^2 \Bigg)^{\frac{1}{2}}\,\text{d}\beta_i \,.
\end{equation}
As both the drift and diffusion functions of the SDE are zero at a consensus state, there is no movement away from consensus or change in the position of consensus once it is established. 

\newpage
\subsection{Summary table}

Table \ref{tab:Summary table} collects the ABMs described in Section \ref{Section: Theory} and their corresponding ODE/SDE limits.
\noindent
\begin{table}[ht!]
    \centering
    \begin{tblr}{
      colspec = {X[1,l]X[3,h]},
      stretch = 0,
      rowsep = 6pt,
      hlines = {black, 1pt},
      vlines = {black, 1pt},
    }
        Agent-based model & Limiting differential equation \\
        Standard ABM (\ref{Subsection: ABM}) & $ \dfrac{dX_i}{dt} = \dfrac{1}{N} \sum\limits_{j=1}^N p_{ij}(X)\, (X_j - X_i)$ \\
        Node degree normalisation (\ref{Subsection: Normalisation}) & $ \dfrac{dX_i}{dt} = \dfrac{1}{k_i} \sum\limits_{j=1}^N A_{ij}\,p_{ij}(X)\, (X_j - X_i)$ \\
        Interaction probability normalisation (\ref{Subsection: Normalisation}) & $  \dfrac{dX_i}{dt} = \Bigg(\sum\limits_{l=1}^N p_{i l}(X)\Bigg)^{-1} \sum\limits_{j=1}^N p_{ij}(X)\, (X_j - X_i)$ \\
        Ambiguity noise (\ref{Subsection: ambiguity noise}) & $ \dfrac{dX_i}{dt} = \dfrac{1}{N} \sum\limits_{j=1}^N p_{ij}(X)\, (X_j - X_i)$ \\
        External noise (\ref{Subsection: External noise}) & $dX_i = \dfrac{1}{N} \sum\limits_{j=1}^N p_{ij}(X)\, (X_j - X_i) \,\text{d}t + \bigg(\dfrac{m_2(\xi)}{N}\bigg)^{\frac{1}{2}} \,\text{d}\beta_i$ \\
        Adaptation noise (\ref{Subsection: External noise}) & $dX_i = \dfrac{1}{N} \sum\limits_{j=1}^N p_{ij}(X)\, (X_j - X_i) \,\text{d}t + \bigg(\dfrac{m_2(\xi)}{N} \,\sum\limits_{j=1}^N p_{ij}(X)\bigg)^{\frac{1}{2}}\,\text{d}\beta_i$ \\
        Random update distance (\ref{Subsection: Random update distance}) & $dX_i = \dfrac{1}{N} \sum\limits_{j=1}^N p_{ij}(X)\, (X_j - X_i) \,\text{d}t + \Bigg(\dfrac{m_2(\nu)}{N^2} \sum\limits_{j=1}^N p_{ij}(X)\, (X_j - X_i)^2 \Bigg)^{\frac{1}{2}}\,\text{d}\beta_i$ \\
    \end{tblr}
    \caption{Summary of variations on the ABM defined in Section \ref{Subsection: ABM} and the corresponding limiting differential equations. $m_2(\cdot)$ is defined in \eqref{Eqn: m_k}.}
    \label{tab:Summary table}
\end{table}

\section{Numerical results} \label{Section: Numerical Results}

In order to demonstrate the convergence results described in Section \ref{Section: Theory} we perform numerical simulations of both the ABM and ODE/SDE limits. To allow for a fair comparison between the different models we use the same setup throughout. 

Let $N = 50$. We select a set of initial conditions $x(0)$ uniformly at random in the interval $[-1,1]$ for each agent, and use these same initial conditions in all simulations. Let $h = 10^{-5}$ and final time $T = 20$. We use the BC interaction function with Gaussian selection noise, given by
\begin{equation} \label{Eqn: numerical examples interaction function}
    p_{ij}(x) = 1 - F_{\mathcal{N}(0,1)}\bigg( \frac{|x_i - x_j| - R}{0.01} \bigg),
\end{equation}
where $F_{\mathcal{N}(0,1)}$ is the cumulative distribution function of a standard normal distribution. Interaction functions of this form can be seen in the left panel of Figure \ref{fig: example 'noisy' interaction functions}.

In this section we will compare timeseries of ABMs and their corresponding limits. For two timeseries $X$ and $Y$ we calculate the error timeseries according to 
\begin{equation} \label{Eqn: Error timeseries}
    \text{Error}(t) = \sum_{i=1}^N \big|X_i(t) - Y_i(t) \big|.
\end{equation}

\subsection{ODE limits} \label{Subsection: Simulated demonstration}

\begin{figure}[ht!]
    \centering
    \includegraphics[width = \linewidth, trim = {2cm 0.5cm 2cm 1cm}, clip]{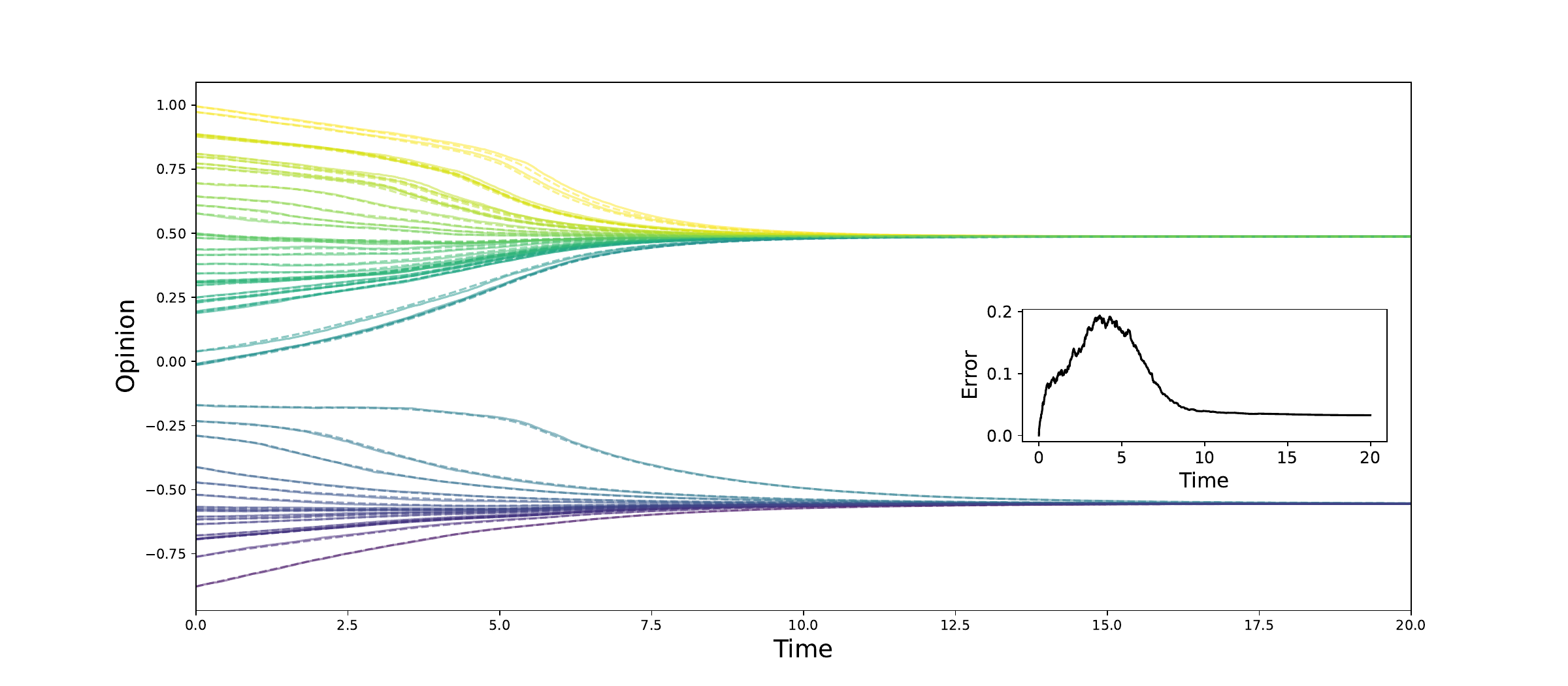}
    \caption{Comparison of the dynamics of the ABM (Section \ref{Subsection: ABM}) and ODE \eqref{Eqn: ODE model} with $h=10^{-5}$. A single realisation of the ABM is plotted in solid lines, with the ODE plotted in dashed lines. Colours represent the initial opinions of the $N=50$ different agents. An inset shows the error, calculated according to \eqref{Eqn: Error timeseries}. There is an excellent match between the two sets of dynamics.}
    \label{fig: example N = 50}
\end{figure}

Solutions to the ODE are calculated numerically using a forward Euler scheme with a timestep of $0.01$. Figure \ref{fig: example N = 50} shows a comparison between a single realisation of the ABM (Section \ref{Subsection: ABM} with update rule \eqref{Eqn: Standard update scheme}) in solid lines and the limiting ODE \eqref{Eqn: ODE model} in dashed lines. Colours represent the initial opinions of the $N=50$ different agents. For this value of $h$ there is clearly an excellent match between the two timeseries, demonstrating the convergence of the ABM to the ODE system. 

To quantify this comparison we examine the difference between the ODE model ABM for various values of $h$. Here the error over the whole timeseries is calculated as follows: realisations of the ABM are sampled at the timesteps of the ODE, the Frobenius norm of the difference between these matrices is then calculated and this value is divided by the length of the simulation to give the error. For each value of $h$ the ABM is simulated 100 times and the error calculated for each simulation. 

\begin{figure}[ht!]
    \centering
    \includegraphics[width = \linewidth, trim = {2cm 0.5cm 2cm 1cm}, clip]{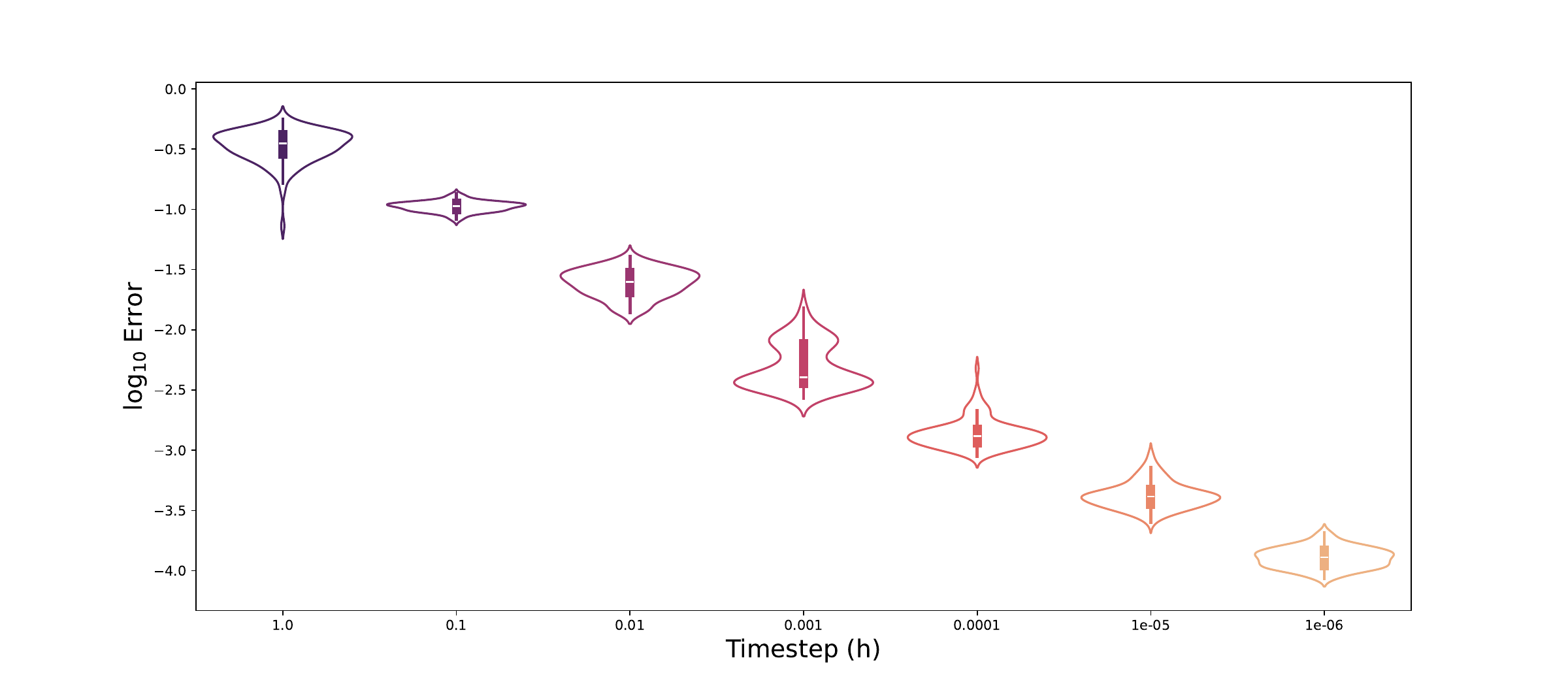}
    \caption{Violin plot of error between the ABM (Section \ref{Subsection: ABM}) and ODE \eqref{Eqn: ODE model} for different values of $h$. For each value of $h$ the ABM is run 100 times with the same initial conditions and the error against the ODE calculated.}
    \label{fig: errors}
\end{figure}

Our results are shown in a violin plot in Figure \ref{fig: errors}. For each value of $h$ the violin outline shows the distribution of error values observed over the 100 realisations, with the mean and interquartile range shown within the violin shape. Figure \ref{fig: errors} shows that as $h$ is decreased, both the mean error and spread in errors decrease. The convergence results we apply in the Supplementary Material do not provide a bound on the convergence rate between the models, only guaranteeing weak convergence in the limit $h\rightarrow0$. However, these results show, as would be expected, that decreasing the value of $h$ gives a better approximation to the ODE limit. 

\subsubsection{Normalisation}

We also implement the two alternative methods for selecting the interacting pair $i$ and $j$ discussed in Section \ref{Subsection: Normalisation}. To consider the limiting ODE \eqref{Eqn: ODE model, k_i normalisation} we must first define a network. We generate an Erd\H{o}s-R\'enyi random network \cite{erdHos1960evolution} and take its adjacency matrix. This adjacency matrix is then adapted by setting $A_{ii}=1$ for all $i=1,\dots,N$. The ABM is simulated again, and now at each timestep individual $i$ is selected first, then $j$ is selected according to \eqref{Eqn: node degree prob}. The results are shown in Figure \ref{fig: example N = 50, network normalisation}. Although the dynamics are more complex than in the case of the fully connected network (as in the original ABM) there is still an excellent match between the ABM and ODE limit \eqref{Eqn: ODE model, k_i normalisation}. 

\begin{figure}[ht!]
    \centering
    \includegraphics[width = \linewidth, trim = {2cm 0.5cm 2cm 1cm}, clip]{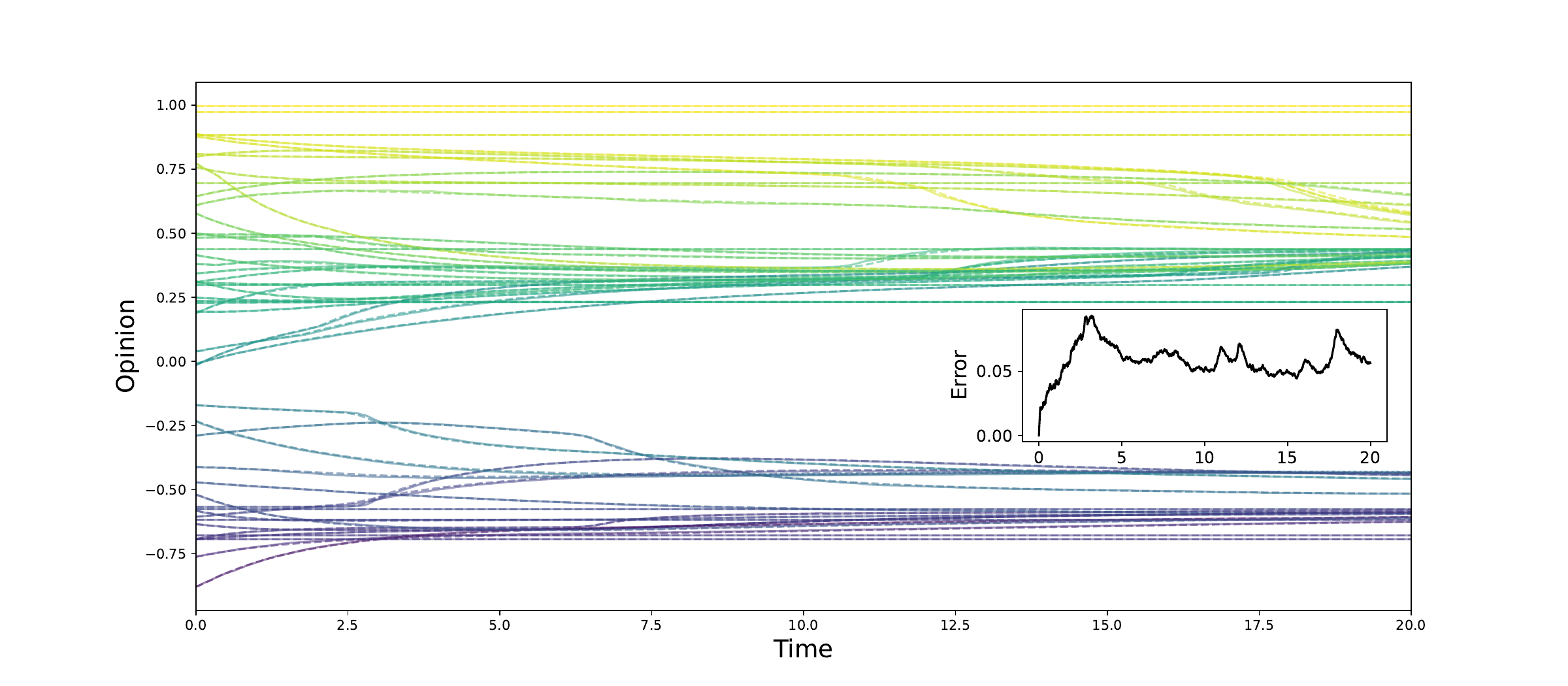}
    \caption{Comparison of the dynamics of the ABM on a network (Section \ref{Subsection: ABM}) in which individual $j$ is chosen according to \eqref{Eqn: node degree prob} and the limiting ODE \eqref{Eqn: ODE model, k_i normalisation} with $h=10^{-5}$. A single realisation of the ABM is plotted in solid lines, with the ODE plotted in dashed lines. Colours represent the initial opinions of the $N=50$ different agents. An inset shows the error, calculated according to \eqref{Eqn: Error timeseries}. There is again an excellent match between the two sets of dynamics.}
    \label{fig: example N = 50, network normalisation}
\end{figure}

Similar results are obtained for the case in which individual $i$ first calculates the probability of interaction with each individual, then selects $j$ with probability proportional to the interaction probability \eqref{Eqn: interaction prob normalisation}. This model does not include a network. These results are shown in Figure \ref{fig: example N = 50, probability normalisation}. Here the population reaches the clustered state much faster than seen in Figure \ref{fig: example N = 50}. As we observe two clusters forming, individuals will have a probability of interaction that is close to $1$ for approximately half the population, and a probability of interaction that is close to $0$ for the rest of the population. Hence the normalisation factor in \eqref{Eqn: ODE model, p_ij normalisation} ($ \sum_{\ell=1}^N p_{i\ell}$) will be approximately half of that in \eqref{Eqn: ODE model} ($N$), so the dynamics are notably faster. This effect is more pronounced at the beginning of the dynamics, as the spread of initial opinions reduces the value of $ \sum_{\ell=1}^N p_{i\ell}$, thus leading to even faster dynamics. It is also more pronounced in individuals with extreme opinions, meaning their movement towards more moderate opinions is even faster. Both the ABM and ODE limit \eqref{Eqn: ODE model, p_ij normalisation} capture these behaviours.  

\begin{figure}[ht!]
    \centering
    \includegraphics[width = \linewidth, trim = {2cm 0.5cm 2cm 1cm}, clip]{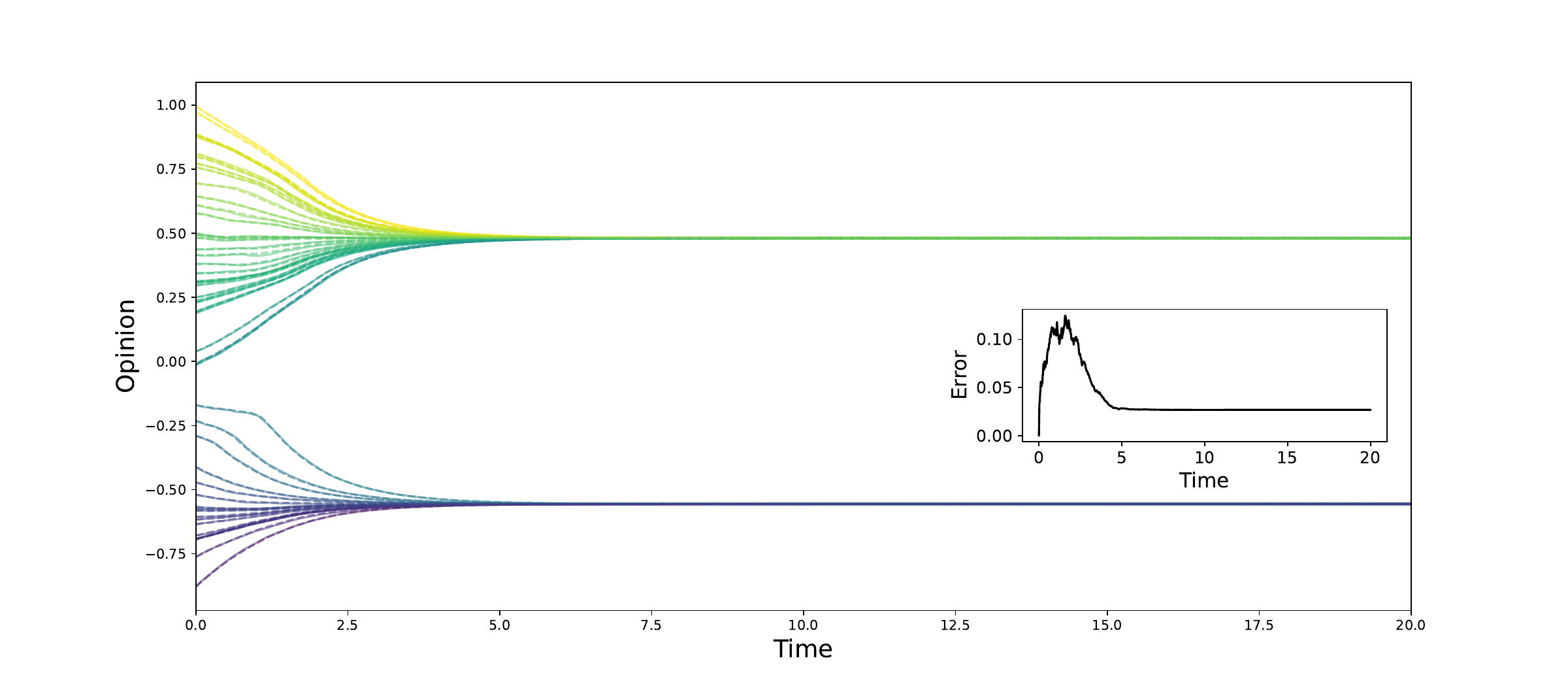}
    \caption{Comparison of the dynamics of the ABM (Section \ref{Subsection: ABM}) in which individual $j$ is chosen according to \eqref{Eqn: interaction prob normalisation} and the limiting ODE \eqref{Eqn: ODE model, p_ij normalisation} with $h=10^{-5}$. A single realisation of the ABM is plotted in solid lines, with the ODE plotted in dashed lines. Colours represent the initial opinions of the $N=50$ different agents. An inset shows the error, calculated according to \eqref{Eqn: Error timeseries}. There is again an excellent match between the two sets of dynamics.}
    \label{fig: example N = 50, probability normalisation}
\end{figure}

We do not include numerical results for the variation in which $i$ and $j$ are chosen without replacement, the variation in which both $i$ and $j$ update their opinions, or the addition of ambiguity noise, as the limiting models are exactly the same as the base ODE case. 

\subsection{SDE limits} \label{Subsection: SDE Numerical examples}

Realisations of the SDE are calculated using an Euler-Maruyama scheme with a timestep of $t = 0.01$. For each of the comparisons in this section we use the same initial conditions and simulation setup as before. We run 5000 realisations and calculate the mean and variance between realisations, over time, for each individual. Due to the multiple sources of stochasticity, this high number of realisations is chosen to allow sufficient opportunity for the variances of each model to be representative of their true values. 

\subsubsection{External noise}

Firstly we consider the ABM with external noise, given by update rule \eqref{Eqn: External noise update scheme}. We take $\xi^h \sim \mathcal{N}(N h, sh)$ with constant $s = 0.05$. Here $\mathds{E}\big[(\nu^h)^2\big] = (N h)^2 + sh$, meaning $m_2(\nu) = s > 0$.  

\begin{figure}[ht!]
    \centering
    \includegraphics[width=\linewidth]{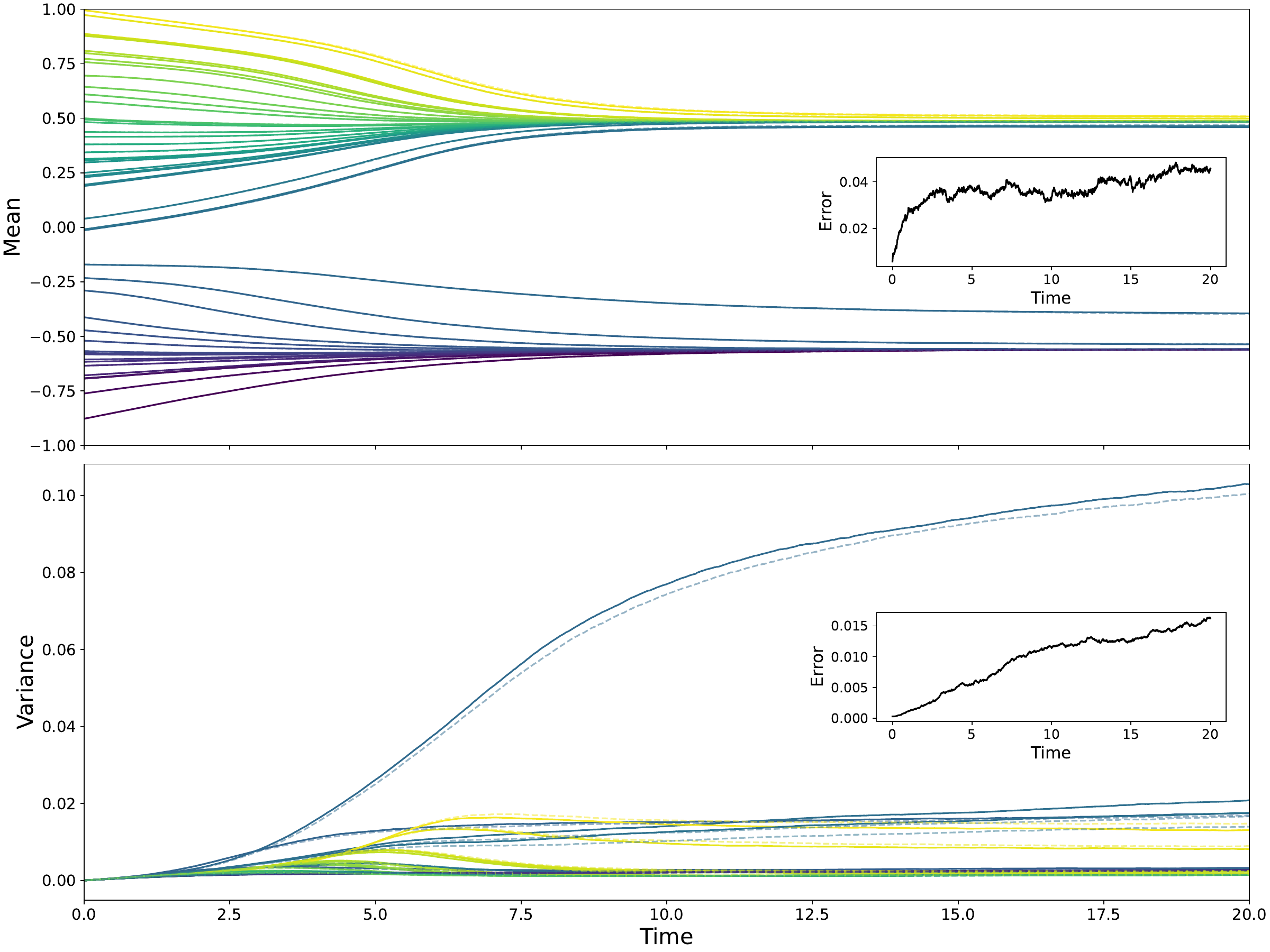}
    \caption{Comparison of the ABM with \textbf{external noise} \eqref{Eqn: External noise update scheme} and limiting SDE system \eqref{Eqn: limiting SDE external noise} with the setup described at the beginning of Section \ref{Subsection: SDE Numerical examples}. Solid lines show the mean (top panel) and variance (bottom panel) in the opinion of each individual across 5000 realisations of ABM, dashed lines show the same for the corresponding SDE system. Both panels include an inset showing the error, calculated according to \eqref{Eqn: Error timeseries}.}
    \label{fig: comparison additive noise}
\end{figure}

\begin{figure}[ht!]

    \begin{subfigure}{\textwidth} 
    \centering
        \begin{subfigure}{.3\textwidth} 
            \centering
            \includegraphics[width=\linewidth, trim = {3cm 2cm 3cm 1cm}, clip]{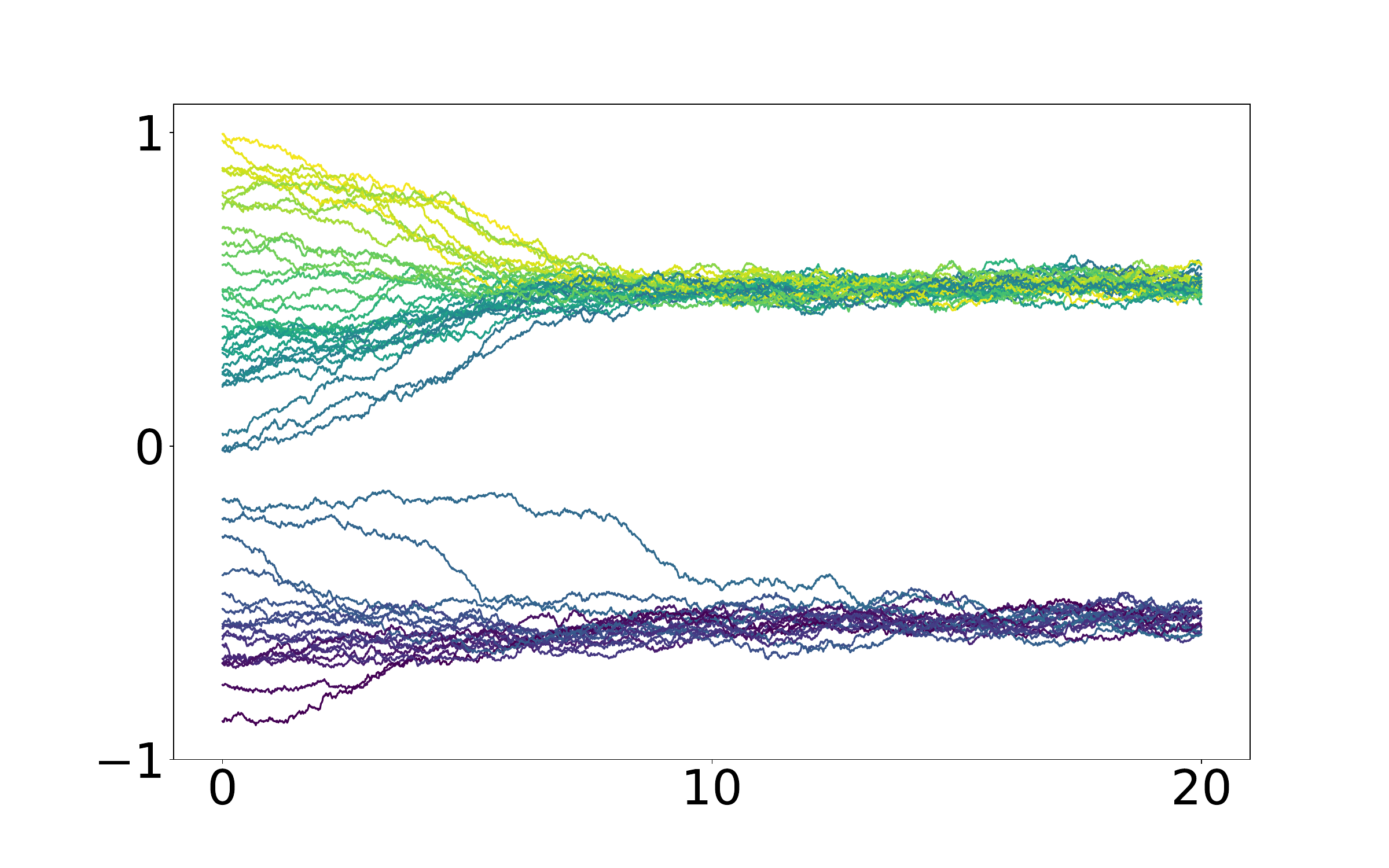}
        \end{subfigure}
        \begin{subfigure}{.3\textwidth} 
            \centering
            \includegraphics[width=\linewidth, trim = {3cm 2cm 3cm 1cm}, clip]{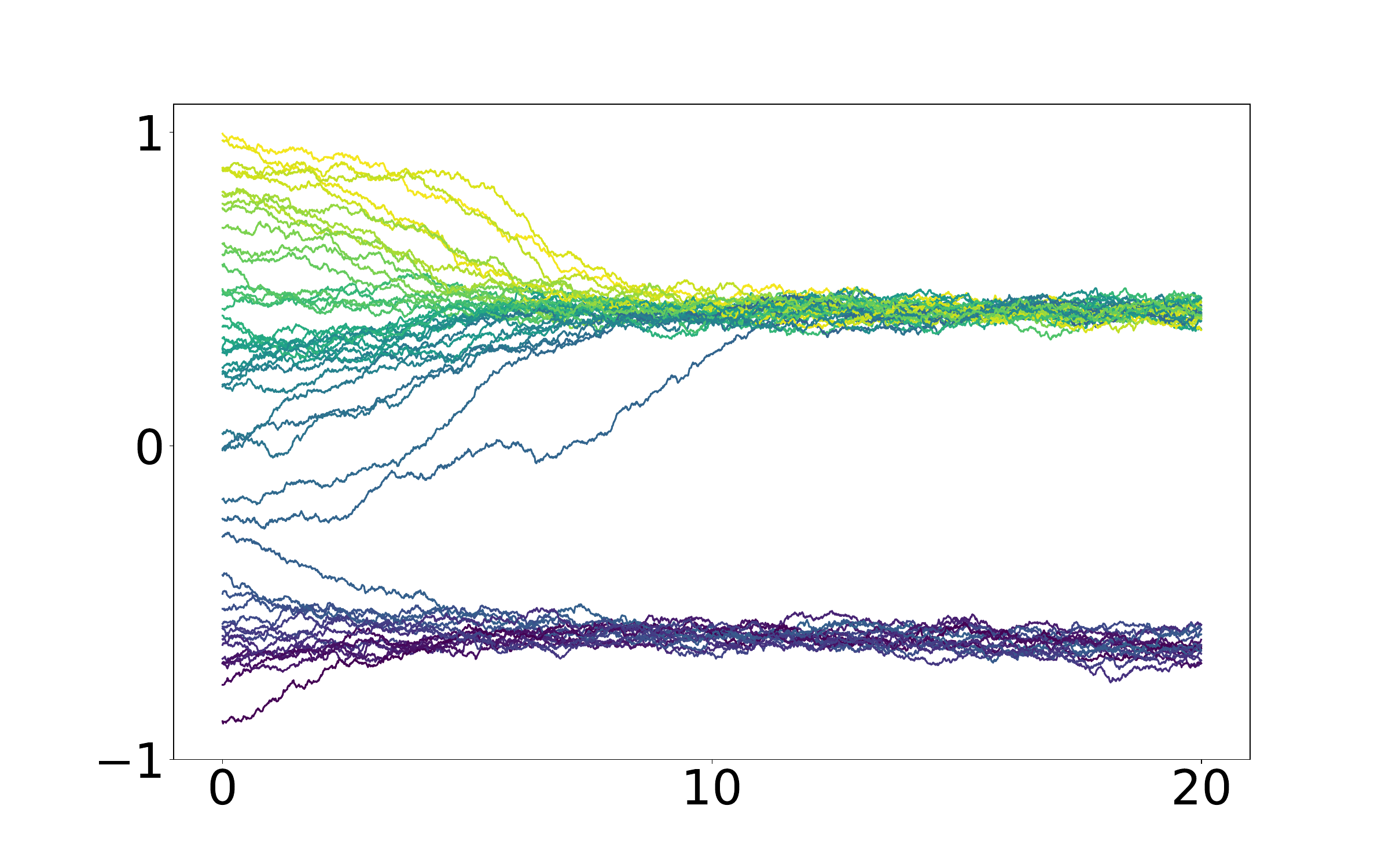}
        \end{subfigure}
        \begin{subfigure}{.3\textwidth} 
            \centering
            \includegraphics[width=\linewidth, trim = {3cm 2cm 3cm 1cm}, clip]{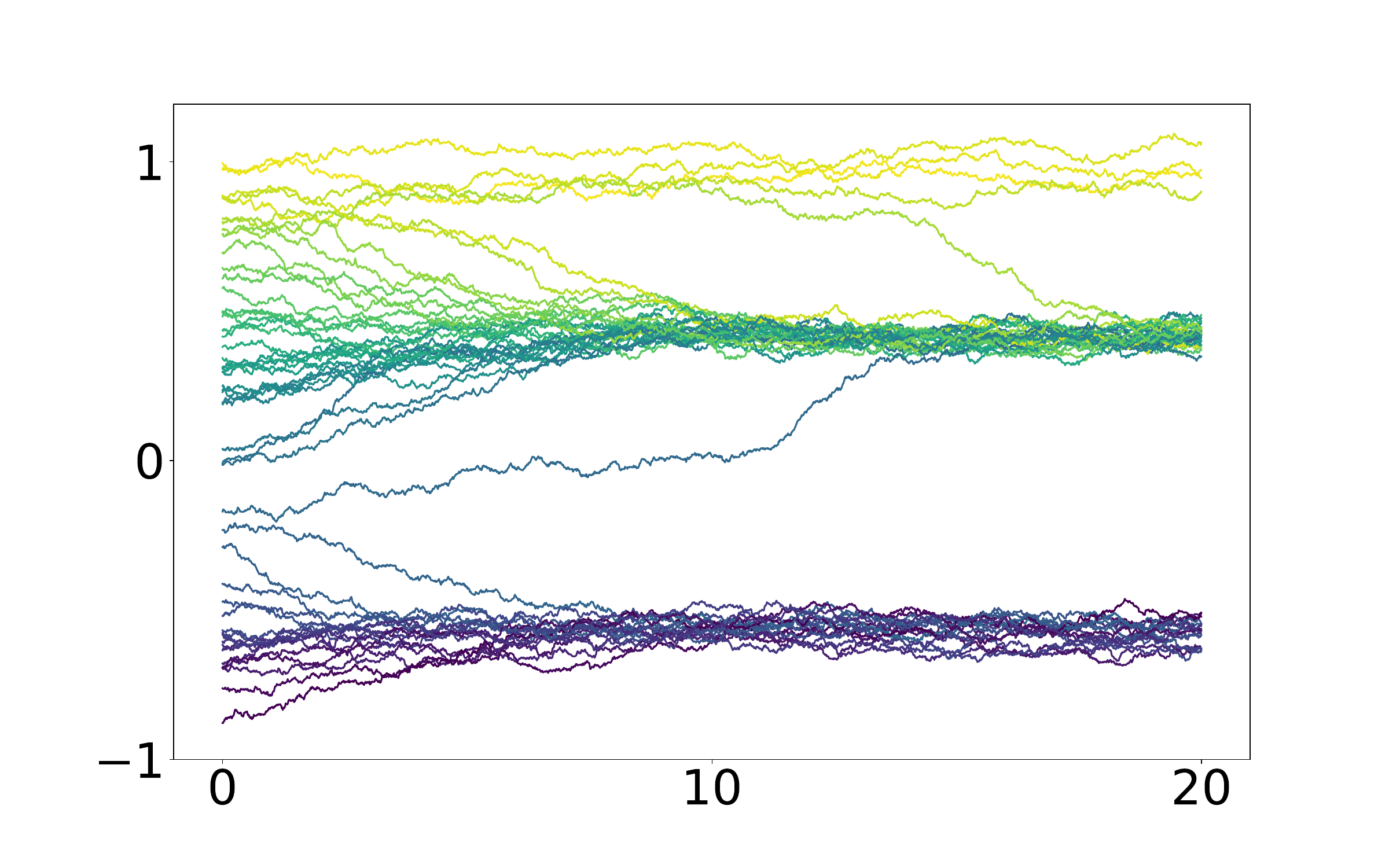}
        \end{subfigure}
        \label{fig: pw examples xi}
    \end{subfigure}

    \caption{Example timeseries of the ABM with \textbf{external noise} \eqref{Eqn: External noise update scheme}, using the setup described at the beginning of Section \ref{Section: Numerical Results}.}
    \label{fig:Example timeseries external}
\end{figure}

We observe excellent agreement between both the means and variances of the ABM and SDE realisations. As in the ODE model, the means reach a clustered state. For the majority of individuals their mean behaviour is very similar to that of the ODE model and the variance between realisations is low. This is perhaps unsurprising given the relatively small size of $m_2(\nu)$. However, for a small number of individuals with initial opinions close to zero (appearing in blue in Figure \ref{fig: comparison additive noise}) their mean behaviour is different and their variance is significantly higher. This indicates that, with the addition of noise, they may sometimes enter either of the two major opinion clusters. There are also some individuals with initial opinions close to one (appearing in yellow in Figure \ref{fig: comparison additive noise}) whose variance is higher, indicating they may sometimes not join the nearby opinion cluster or make take a long time to do so. Both of these behaviours can be observed in the selected example timeseries in Figure \ref{fig:Example timeseries external}.

These example timeseries also show that, unlike the deterministic case, individuals' opinions continue to change after the opinion clusters have formed. This is due to the fact that interactions are still occurring, meaning individuals' opinions are still affected by the external noise, even if their opinions are very close to those they are interacting with. This can also be seen in the variances in the lower panel of Figure \ref{fig: comparison additive noise}, which continue to increase as the simulation goes on. If these simulations were allowed to run indefinitely these variances would continue to increase as the opinion clusters followed a random walk, driven by the continued external noise. 

\subsubsection{Adaptation noise}

We now consider the ABM with adaptation noise, given by the update rule \eqref{Eqn: External noise with interactions update scheme}. As previously we take $\xi^h \sim \mathds{N}(N h, sh)$ with constant $s = 0.05$. 

\begin{figure}[ht!]
    \centering
    \includegraphics[width=\linewidth]{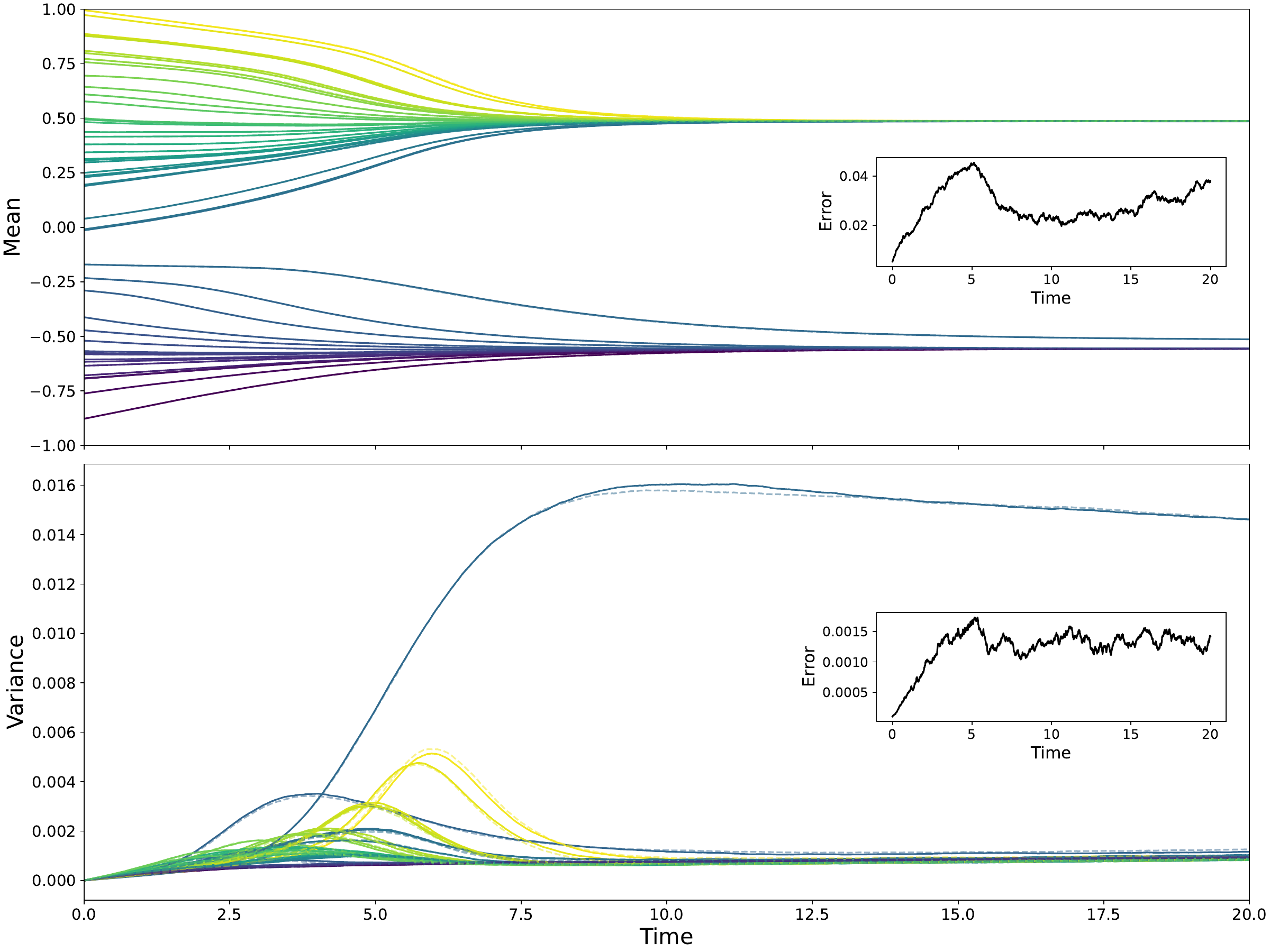}
    \caption{Comparison of the ABM with \textbf{adaptation noise} \eqref{Eqn: External noise with interactions update scheme} and limiting SDE system \eqref{Eqn: limiting SDE external noise on interactions} with the setup described at the beginning of Section \ref{Subsection: SDE Numerical examples}. Solid lines show the mean (top panel) and variance (bottom panel) in the opinion of each individual across 5000 realisations of ABM, dashed lines show the same for the corresponding SDE system. Both panels include an inset showing the error, calculated according to \eqref{Eqn: Error timeseries}.}
    \label{fig: comparison ambiguity noise}
\end{figure}

In this case we again observe (in Figure \ref{fig: comparison ambiguity noise}) an excellent match between the behaviours of the ABM and limiting SDE system \eqref{Eqn: limiting SDE external noise on interactions}. The behaviour of the means (across realisations) is again similar to that of the ODE model in Figure \ref{fig: example N = 50}. As in the previous case (external noise) we observe one individual with significantly higher variance than all others, again explained by entering different opinion clusters, or neither opinion cluster, in different realisations. In this case the variance of this individual drops slightly towards the end of the simulations, indicating that they do eventually join one of the clusters. 

\begin{figure}[ht!]

    \begin{subfigure}{\textwidth} 
    \centering
        \begin{subfigure}{.3\textwidth} 
            \centering
            \includegraphics[width=\linewidth, trim = {3cm 2cm 3cm 1cm}, clip]{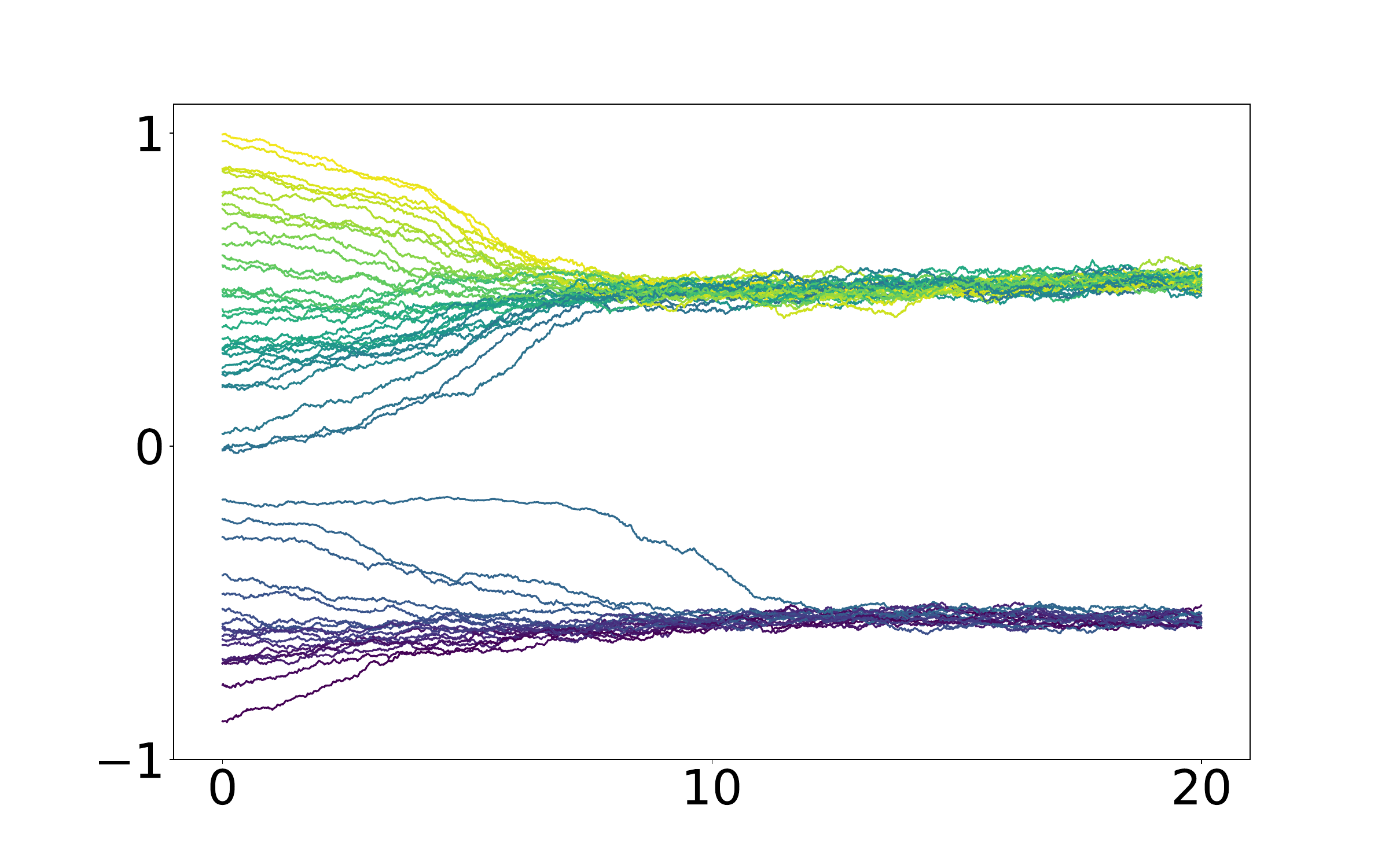}
        \end{subfigure}
        \begin{subfigure}{.3\textwidth} 
            \centering
            \includegraphics[width=\linewidth, trim = {3cm 2cm 3cm 1cm}, clip]{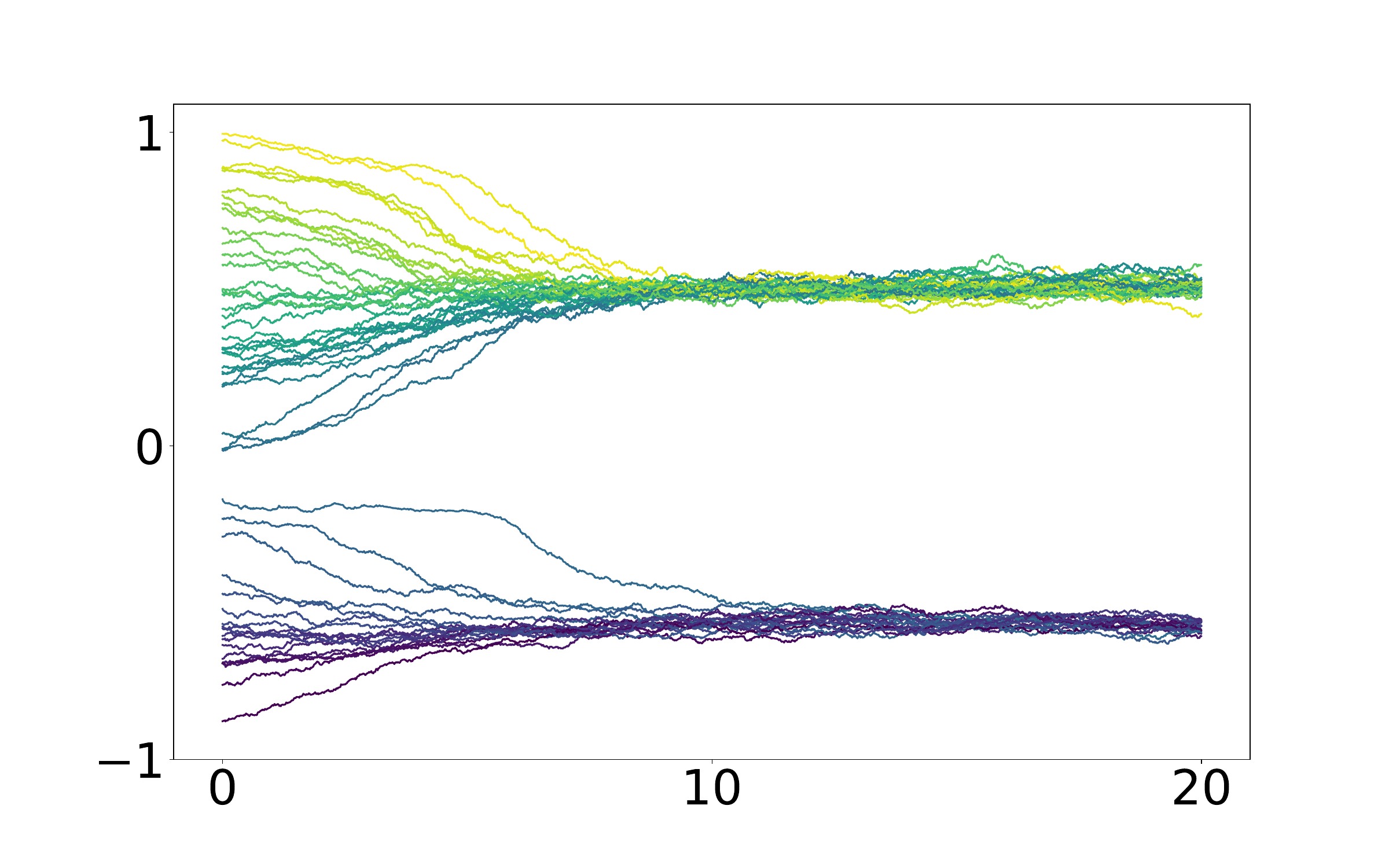}
        \end{subfigure}
        \begin{subfigure}{.3\textwidth} 
            \centering
            \includegraphics[width=\linewidth, trim = {3cm 2cm 3cm 1cm}, clip]{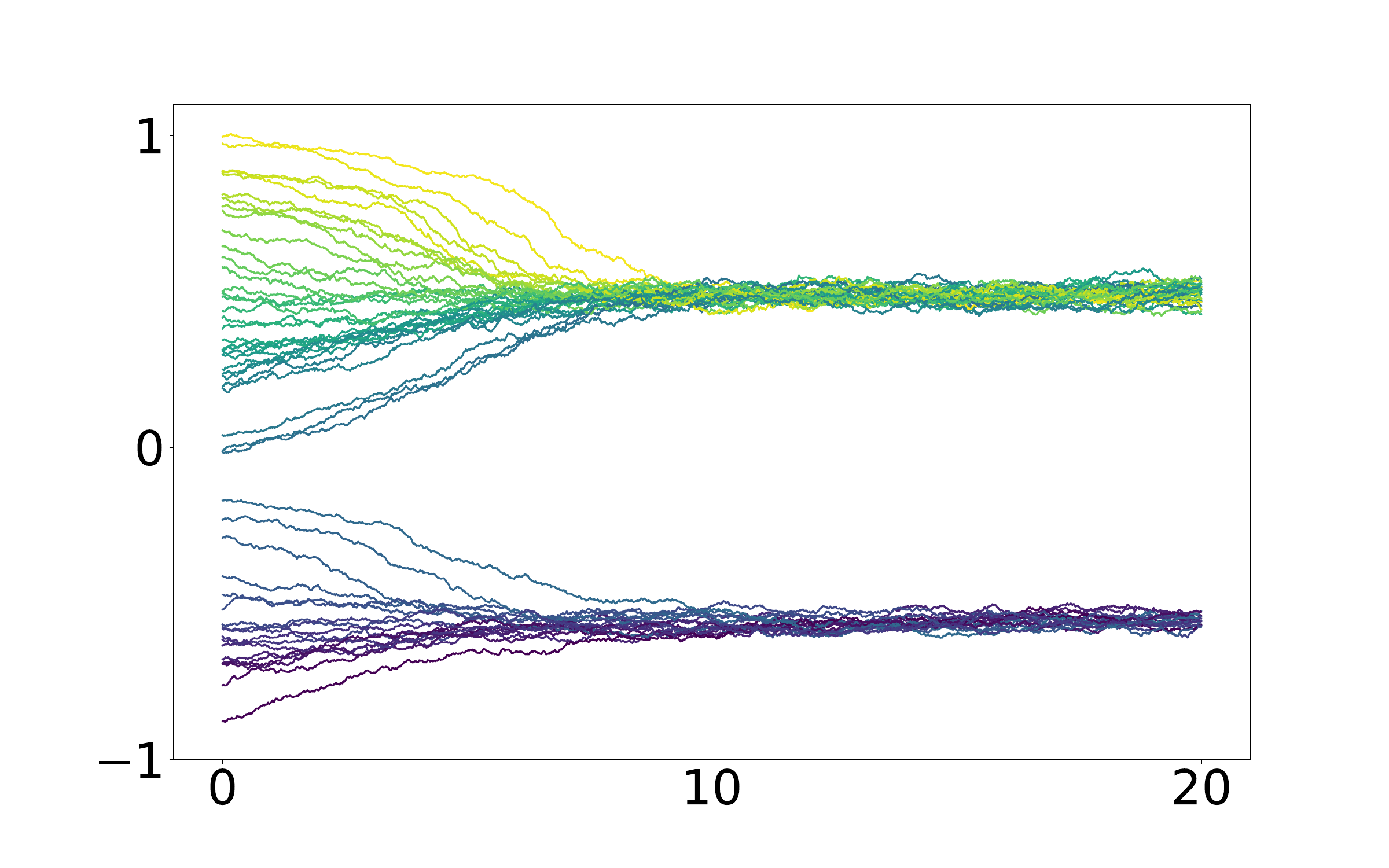}
        \end{subfigure}
        
        \label{fig: pw examples adaptation}
    \end{subfigure}
    
    \caption{Example timeseries of the ABM with \textbf{adaptation noise} \eqref{Eqn: External noise with interactions update scheme}, using the setup described at the beginning of Section \ref{Section: Numerical Results}.}
    \label{fig:Example timeseries adaptation}
\end{figure}

A similar behaviour can be observed in the individuals with initial opinion near one (appearing in yellow) whose variance initially grows, peaks around $t = 6$, then drops again. This indicates that these individuals may take some time to join an opinion cluster, but always do so eventually. This variation in the time taken for clusters to form can be seen in the example timeseries in Figure \ref{fig:Example timeseries adaptation}. As the population splits into two non-interacting clusters, and individuals only experience noise when they interact, this model has comparatively smaller noise than in the previous case, leading to more predictable outcomes. 

\subsubsection{Random update distance}

Here we give a numerical demonstration of the convergence to the SDE system \eqref{Eqn: limiting SDE noisy update distance} shown in Section \ref{Subsection: Random update distance}. Let $\nu^h \sim \mathds{N}(N h, sh)$, with constant $s = 5$. As the noise term $\nu$ is multiplied by the distance between individuals' opinions, we choose a larger value of $s$ so that the presence of noise is still clear.

\begin{figure}[ht!]
    \centering
    \includegraphics[width=\linewidth]{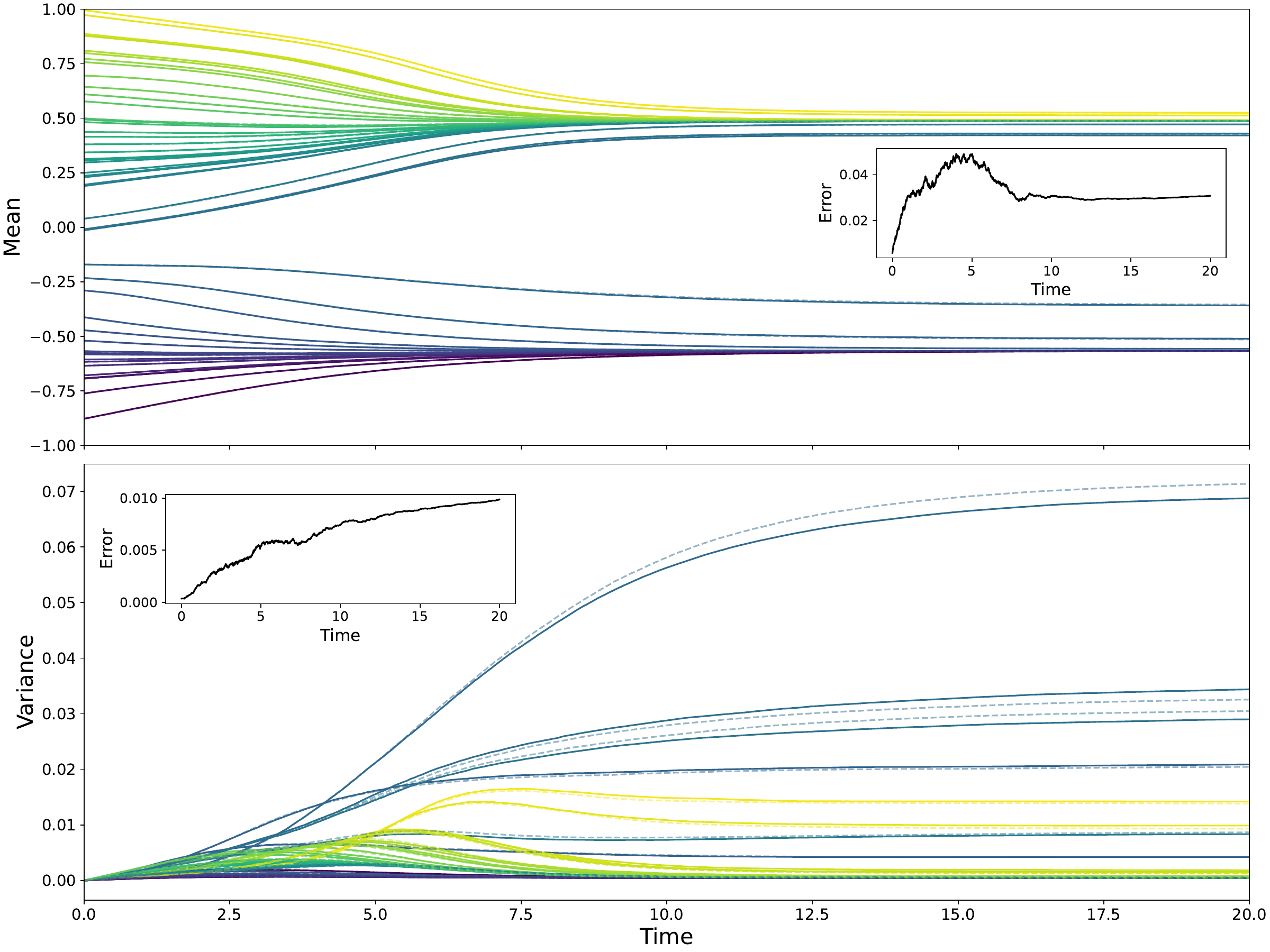}
    \caption{Comparison of the ABM with \textbf{random update distance} \eqref{Eqn: Noisy update distance update scheme} and limiting SDE system \eqref{Eqn: limiting SDE noisy update distance} with the setup described at the beginning of Section \ref{Subsection: SDE Numerical examples}. Solid lines show the mean (top panel) and variance (bottom panel) in the opinion of each individual across 5000 realisations of ABM, dashed lines show the same for the corresponding SDE system. Both panels include an inset showing the error, calculated according to \eqref{Eqn: Error timeseries}.}
    \label{fig: comparison noisy update distance}
\end{figure}

The top panel of Figure \ref{fig: comparison noisy update distance} compares the mean opinion of each agent in the ABM in solid lines and the SDE system in dashed lines. There is clearly a good fit between the two. Similarly the behaviours of the variances shown in the lower panel of Figure \ref{fig: comparison noisy update distance} match well, although the values towards the end of the timeseries differ slightly for some individuals. 

The behaviour of the means and variances mirrors that seen in Figure \ref{fig: comparison additive noise} for the ABM and SDE system with external noise. However, in this case, there are more individuals with initial opinions near zero who have a high variance. This may be caused by the possibility of $\nu$ taking negative values, meaning interactions can cause individuals' opinions to move apart, adding to the uncertainty on which opinion cluster individuals will enter. 

\begin{figure}[ht!]

    \begin{subfigure}{\textwidth} 
    \centering
        \begin{subfigure}{.3\textwidth} 
            \centering
            \includegraphics[width=\linewidth, trim = {3cm 2cm 3cm 1cm}, clip]{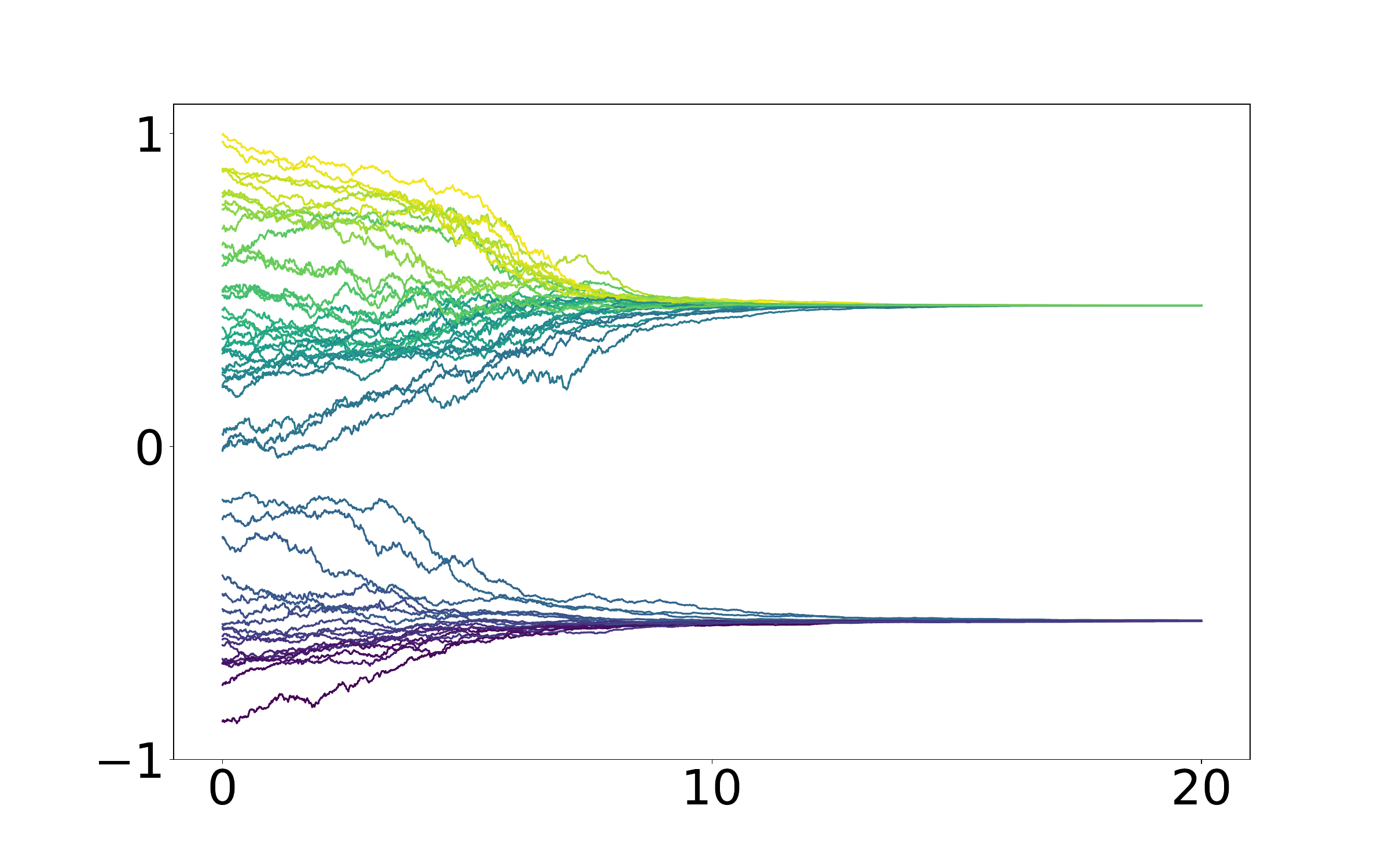}
        \end{subfigure}
        \begin{subfigure}{.3\textwidth} 
            \centering
            \includegraphics[width=\linewidth, trim = {3cm 2cm 3cm 1cm}, clip]{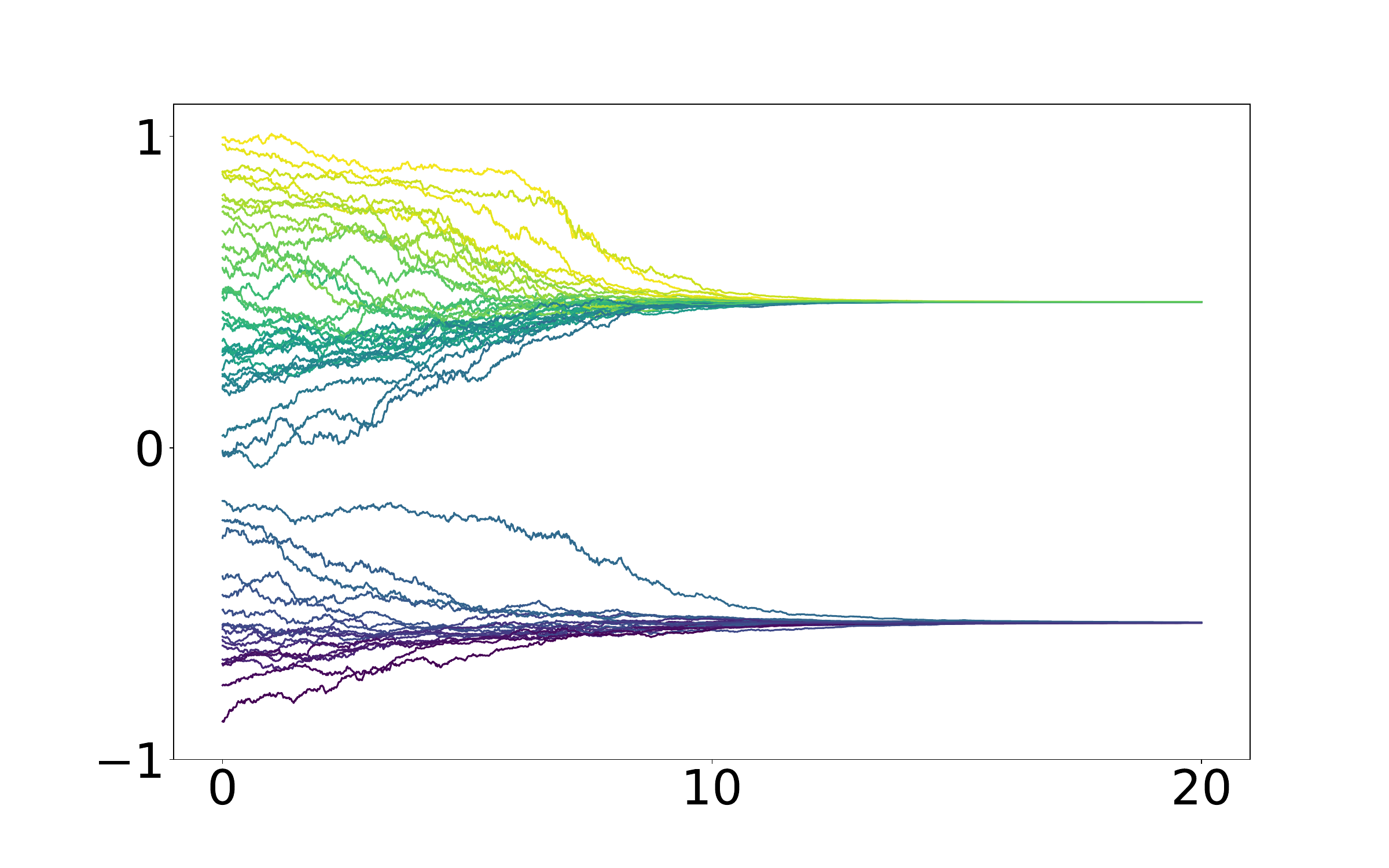}
        \end{subfigure}
        \begin{subfigure}{.3\textwidth} 
            \centering
            \includegraphics[width=\linewidth, trim = {3cm 2cm 3cm 1cm}, clip]{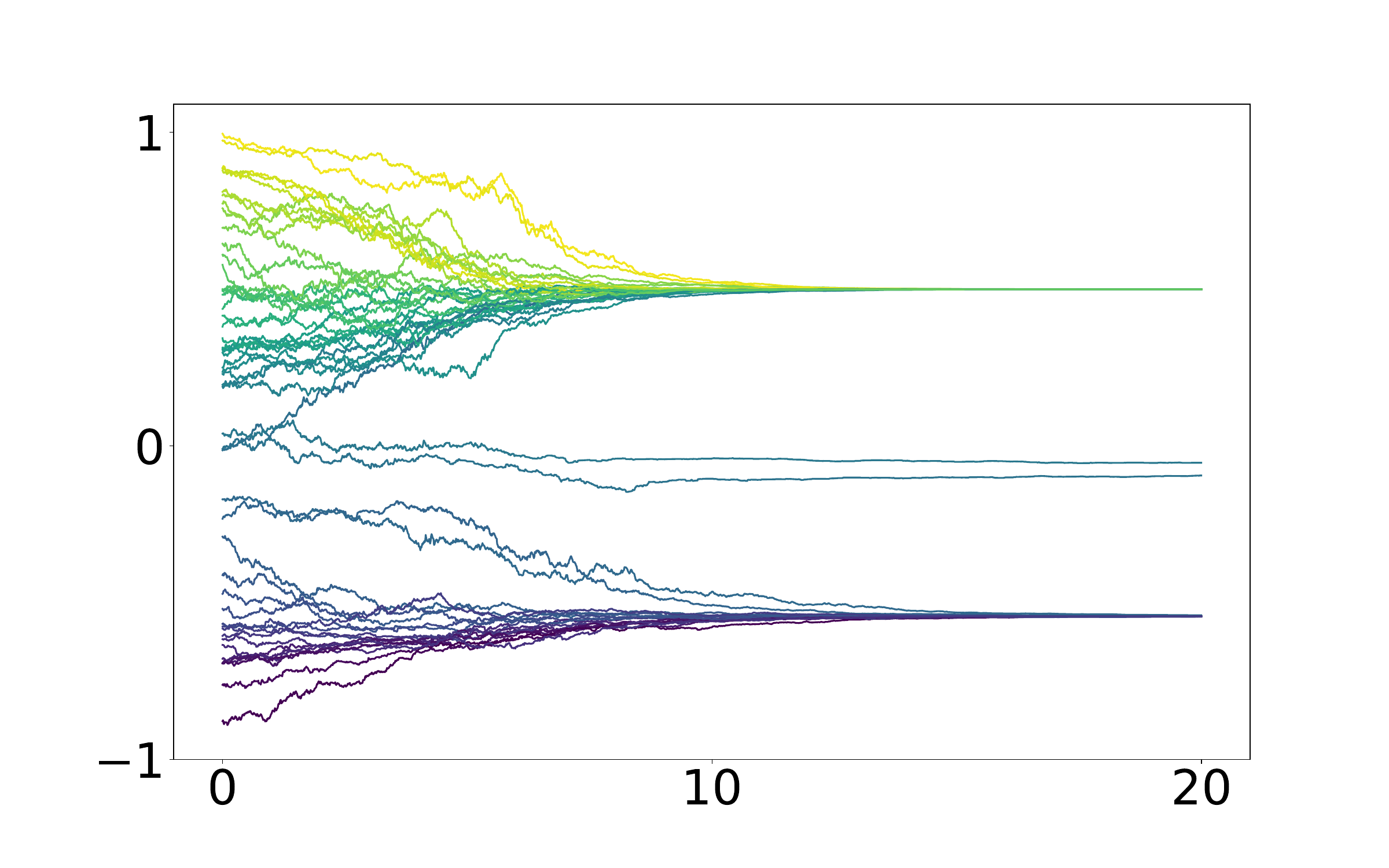}
        \end{subfigure}

        \label{fig: pw examples}
    \end{subfigure}

    \caption{Example timeseries of the ABM with \textbf{random update distance} \eqref{Eqn: Noisy update distance update scheme}, using the setup described at the beginning of Section \ref{Section: Numerical Results}.}
    \label{fig:Example timeseries nu}
\end{figure}

There is a key difference in the example timeseries shown in Figure \ref{fig:Example timeseries nu}, compared against those in Figures \ref{fig:Example timeseries external} and \ref{fig:Example timeseries adaptation}: towards the end of each simulation, individuals' opinions stop changing. This is due to the fact that the noise term $\nu$ is multiplied by the distance $(x_j - x_i)$. Towards the end of the simulation, pairs of individuals either have very similar opinions, meaning $(x_j - x_i)$ is very small, or no longer interact. Hence the noise arising from $\nu$ also becomes very small and the variance plateaus. In this case, if the simulations were left to run indefinitely there would be no further change in individuals' opinions and no continued increase in variance, as once a clustered state is reached it is fixed. This also makes the dynamics less predictable, as there is no possibility that individuals will eventually, due to continued stochasticity, join a nearby opinion cluster. 

\section{Discussion} \label{Section: Discussion}

The theoretical results presented in Section \ref{Section: Theory} provide a framework in which choices made in the ABM setting can be compared directly to those made in the ODE/SDE setting. We show that under a simultaneous rescaling of time and update distance the behaviour of the ABM mirrors that of DEMs. This confirms the assumptions inherent to DEMs: that pairwise interactions between individuals produce very small changes in opinion, but occur very frequently. These assumptions are not always valid in real-world scenarios, but when they are we have shown that DEMs are a very good approximation of ABMs. Moreover, we show how including external noise or a random update distance can translate into stochasticity in the limiting model, motivating the choice of diffusion in an SDE model. Our numerical results in Section \ref{Section: Numerical Results} demonstrate the excellent match between the behaviour of the various ABM and SDE models, showing how both can capture a variety of complex behaviours. This means that the benefits of DEMs, such as reduced computational cost and increased analytic tractability, can be combined with the interpretability of ABMs. 

However, the methods we have presented are not without their limitations. As previously discussed, we cannot take $p_{ij}$ to be the BC interaction function without selection noise, as this discontinuity prevents the necessary guarantee of existence and uniqueness of solutions to the corresponding DEM. This issue occurs whenever we wish to implement the BC interaction function in continuous time \cite{ceragioli2021generalized} and questions of existence and uniqueness of different types of solutions to the ODE model have been addressed \cite{ceragioli2012continuous,blondel2010continuous}, but further technical work is required to allow for a comparison against the ABM. In addition, this approach cannot capture effects such as Pineda noise \cite{pineda2009noisy} in which jumps in opinions cannot be bounded by a function of $h$, although convergence to a jump-diffusion process could be considered. 

Figure \ref{fig:Example timeseries external} shows an example timeseries in which noise causes some individuals' opinions to leave the initial interval $[-1,1]$. Boundary conditions could be defined for the ABM and the changes tracked through to identify when a limiting SDE exists. As there are many different ways that boundary conditions could be enforced (for example: reflecting or absorbing boundary conditions at the ends of the opinion interval; restricting the random variables $\eta$, $\xi$ and $\nu$ to prevent movements outside the convex hull; multiplying noise terms by a function of $x_i$ that is zero at the boundary) we present instead more general results that could be adapted further if necessary.   

This paper confirms the strong connection between agent-based models and differential equations models. It provides the so-far missing mathematical underpinning, which will lead to an increased understanding of these two very popular modeling approaches and an improved communication among communities.

\bibliography{bibliography.bib}
\bibliographystyle{unsrt}

\section*{Acknowledgements}

AN was supported by the Engineering and Physical Sciences Research Council through the Mathematics of Systems II Centre for Doctoral Training at the University of Warwick (reference EP/S022244/1). MTW is partly supported by the Royal Society International Exchange Grant IES/R3/213113.

The authors would also like to thank Peter Steiglechner and Guillaume Deffuant, whose presentations at ODCD23 provided the initial motivation for this work. 

For the purpose of open access, the authors have applied a Creative Commons Attribution (CC-BY) license to Any Author Accepted Manuscript version arising from this submission.
\newpage
\appendix
\section{Supplementary Material} \label{Supplementary Material}

\subsection{Results from Durrett's Stochastic Calculus}

We begin by formulating the ABM as a Markov process. For each $h>0$ let $\Pi^h(x,\cdot)$ be the transition function on $\mathds{R}^N$ and $Y^h_{mh},\,m=0,1,2,\dots$ the corresponding discrete time Markov process. Define $X^h_t = Y^h_{h\lfloor t/h \rfloor}$, i.e. $X^h_t$ is constant on intervals $[mh,(m+1)h]$, as well as 
\begin{align}
    a^h_{ij}(x) &= \frac{1}{h} \int_{|y-x|\leq1} (y_i - x_i)(y_j - x_j)\,\Pi^h(x,dy)\,, \label{eqn: a^h} \\[0.5em]
    b^h_{i}(x) &= \frac{1}{h} \int_{|y-x|\leq1} (y_i - x_i)\,\Pi^h(x,dy)\,, \label{eqn: b^h}
\end{align}
and 
\begin{equation} \label{eqn: jump prob}
    \Delta^h_\varepsilon(x) = \frac{1}{h} \Pi^h\big(x,\mathds{R}^d\backslash B(x,\varepsilon)\big).
\end{equation}
In addition, define functions $a:\mathds{R}^d \rightarrow S_N$ and $b:\mathds{R}^N \rightarrow \mathds{R}^N$, where $S_N$ is the set of non-negative definite real $N\times N$ matrices. 

\begin{assumption} \label{Assumption group: a and b} The functions $a$ and $b$ are continuous and the martingale problem with coefficients $(a,b)$ is well-posed.
\end{assumption}

Assumption \ref{Assumption group: a and b} is satisfied if $a$ and $b$ are both globally Lipschitz continuous, see \cite{durrett2018stochastic}. The stochastic process $X$ being a solution to the martingale problem for coefficients $a$ and $b$ is equivalent to it being a weak solution of the corresponding SDE
\begin{equation}
    dX_t = b(X_t) \,dt + \sqrt{a(X_t)}\,d\beta_t.
\end{equation}
As we are concerned with weak convergence of the ABM to the SDE, we use this representation. See Section 5.4, specifically the discussion of Theorem 4.5, in \cite{durrett2018stochastic} for full details.

\begin{assumption} \label{Assumption group: thm 7.1} For each $i,j$ with $1 \leq i,j \leq N$, $R< \infty$ and $\varepsilon > 0$
\begin{enumerate}[label=\alph*)]
    \item $\lim\limits_{h\searrow0} \sup\limits_{|x|\leq R} | a_{ij}^h(x) - a_{ij}(x) | = 0$ ,\label{Assumption: a^h to a}
    \item $\lim\limits_{h\searrow0} \sup\limits_{|x|\leq R} | b_{i}^h(x) - b_{i}(x) | = 0$ ,\label{Assumption: b^h to b}
    \item $\lim\limits_{h\searrow0} \sup\limits_{|x|\leq R} \Delta^h_\varepsilon(x) = 0$ .\label{Assumption: no jumps}
\end{enumerate}
\end{assumption}

Before stating the main theorem we give definitions of two relevant notions of convergence for stochastic processes. We will consider processes with trajectories in the space of c\`adl\`ag functions, $D([0,T],\mathds{R}^N)$ for any $T < \infty$. Denote by $d$ the metric on this space, as described in \cite{durrett2018stochastic,billingsley2013convergence} (see these for a complete description of this space and its properties). The following definitions of convergence are adapted from those given in \cite{billingsley2013convergence}.

\begin{definition}
    For a random element $X_t$ and a sequence of random elements $\{X^h_t\}$, let $P$ and $P^h$ denote their respective distributions. We say $\{X^h_t\}$ \textbf{converges in distribution} to $X_t$ if $P^h f \rightarrow P f$ for each bounded, continuous real-valued function $f$ on $D([0,T],\mathds{R}^N)$. This is denoted $X^h_t \Rightarrow X_t$.
\end{definition}

\begin{definition}
    For $c\in D([0,T],\mathds{R}^N)$, we say $\{X^h_t\}$ \textbf{converges in probability} to $c$ if, as $h\rightarrow0$,
    \begin{equation}
        \mathds{P}[ d(X^h_t,c) < \epsilon ] \rightarrow 1.
    \end{equation}
\end{definition}
For a given (deterministic) $c\in D([0,T],\mathds{R}^N)$, convergence in probability is equivalent to convergence in distribution \cite{billingsley2013convergence}.

We now present the main results from \cite{durrett2018stochastic} we will use. 

\textbf{Theorem 7.1 statement:} If $X_0^h = x^h \rightarrow x$ and Assumptions \ref{Assumption group: a and b} and \ref{Assumption group: thm 7.1} hold, we have $X_t^h \Rightarrow X_t$ as $h\rightarrow0$, where $X_t$ is the solution of the martingale problem with coefficients $a$ and $b$, with $X_0 = x$.

\subsubsection{Useful Lemmas}

The first lemma replaces the truncated moments in \eqref{eqn: a^h} and \eqref{eqn: b^h} with regular moments and replaces \eqref{eqn: jump prob} with an easier expression to calculate. Let $1\leq i,j \leq N$ and $p \in \mathds{N}$, define
\begin{align}
    \hat{a}^h_{ij}(x) &= \frac{1}{h} \int_{\mathds{R}^N} (y_i - x_i)(y_j - x_j)\,\Pi^h(x,dy)\,, \label{eqn: a^h 2} \\[0.5em]
    \hat{b}^h_{i}(x) &= \frac{1}{h} \int_{\mathds{R}^N} (y_i - x_i)\,\Pi^h(x,dy)\,, \label{eqn: b^h 2}
\end{align}
and 
\begin{equation} \label{eqn: jump prob 2}
    \gamma^h_p(x) = \frac{1}{h} \int_{\mathds{R}^N} |y_i - x_i|^p\,\Pi^h(x,dy)\,.
\end{equation}

\begin{remark}
    As we never apply Theorem 7.1 directly and never calculate $a_{ij}^h$ and $b_{i}^h$, only $\hat{a}^h_{ij}$ and $\hat{b}^h_{i}$, we omit the 'superscript hats' from now on.
\end{remark}

\textbf{Lemma 8.2}: If $p\geq 2$ and for all $0 \leq i,j \leq N$,  $R<\infty$, we have
\begin{enumerate}[label=\alph*)]
    \item $\lim\limits_{h\searrow0} \sup\limits_{|x|\leq R} | a_{ij}^h(x) - a_{ij}(x) | = 0$ ,
    \item $\lim\limits_{h\searrow0} \sup\limits_{|x|\leq R} | b_{i}^h(x) - b_{i}(x) | = 0$ ,
    \item $\lim\limits_{h\searrow0} \sup\limits_{|x|\leq R} \gamma^h_p(x) = 0$ ,
\end{enumerate}
then Assumption \ref{Assumption group: thm 7.1} hold.

A second lemma is helpful in the case of deterministic limits. 

\textbf{Lemma 8.5}: If for all $0 \leq i,j \leq N$, $R<\infty$ we have
\begin{enumerate}[label=\alph*)]
    \item $\lim\limits_{h\searrow0} \sup\limits_{|x|\leq R} | a_{ij}^h(x) | = 0$ , 
    \item $\lim\limits_{h\searrow0} \sup\limits_{|x|\leq R} | b_{i}^h(x) - b_{i}(x) | = 0$ ,
\end{enumerate}
then Assumption \ref{Assumption group: thm 7.1} hold with $a_{ij}(x)\equiv0$.

In this case, as the limit $X_t$ is deterministic, the process $X^h_t$ converges to $X_t$ both in distribution and in probability. 

\subsection{Convergence results} \label{Appendix: Convergence of ABM}

Consider the ABM defined in Section \ref{Subsection: ABM}. We will use Lemma 8.5 from \cite{durrett2018stochastic} (given above) to show that this Markov process converges weakly to the ODE \eqref{Eqn: ODE model}. In the following calculations we assume that all subscript indices $i,j$ satisfy $1 \leq i,j \leq N$.

Let $x$ be the current opinions of the population. As there are $N^2$ possible interactions, there are also $N^2$ states to which the system could move. These are given by
\begin{equation}
    x' = x + e_i \mu^h (x_j - x_i) \,,
\end{equation}
for $i,j=1,\dots,N$, where $e_i$ is the $i^\text{th}$ standard basis vector of $\mathds{R}^N$. This state is obtained if $i$ and $j$ are chosen to interact, which occurs with probability $N^{-2}$, and indeed go on to interact, which occurs with probability $p_{ij}(x)$. If the same individual is chosen as both $i$ and $j$, which occurs with probability $N^{-1}$, the system will not change state. Additionally, if $i\neq j$ are chosen but do not interact, which occurs with probability $1 - p_{ij}(x)$, the system will also not change state. Hence the transition probability function $\Pi^h$ is given by 
\begin{equation} \label{Eqn: Pi^h}
    \Pi^h(x,y) = 
    \begin{cases}
        \dfrac{1}{N^2} p_{ij}(x) & \text{if } y = x +  e_i \mu^h (x_j - x_i) \text{ for some $i\neq j$}, \\[1em]
        \dfrac{1}{N} + \dfrac{1}{N^2}  \sum\limits_{i\neq j} \big(1 - p_{ij}(x)\big) & \text{if } y = x,  \\[1.2em]
        0 & \text{otherwise.}  
    \end{cases}
\end{equation}
Note that $\Pi^h$ depends on $h$ through the update distance $\mu^h$, which determines which states can be reached in one update. 

In order to apply Lemma 8.5 we calculate $a^h_{ij}$ and $b^h_{i}$. These can be understood as approximations to the coefficients $a_{ii}$ and $b_i$ of a limiting SDE model that we wish to determine. If $i\neq j$ then there are no interactions for which both $(y_i - x_i)$ and $(y_j - x_j)$ are both nonzero, since the opinion of only one individual is changed at each timestep. Hence $a^h_{ij}\equiv0$ if $i\neq j$. For $i=1,\dots,N$, we calculate the following integrals by summing over the possible interactions, noting that only those interactions involving individual $i$ appear as we multiply by the change in the $e_i$ direction. Recall that we set $\mu = Nh$, and therefore
\begin{align*}
    b^h_{i}(x) 
    &= \frac{1}{h} \sum_{j=1}^N \mu^h (x_j - x_i)\,\frac{1}{N^2} p_{ij}(x)\,,  
    = \frac{1}{N} \sum_{j=1}^N p_{ij}(x)\, (x_j - x_i)\, ,  \\
    a^h_{ii}(x) 
    &= \frac{1}{h} \sum_{j=1}^N \big(\mu^h (x_j - x_i)\big)^2\,\frac{1}{N^2} p_{ij}(x)\,,  
    = h \Bigg( \sum_{j=1}^N p_{ij}(x)\,(x_j - x_i)^2\, \Bigg) \,, 
\end{align*}
Using these we now check the conditions of Lemma 8.5. Since $b_i^h = b_i$, as seen in the limiting ODE \eqref{Eqn: ODE model}, the second is immediately satisfied. Hence it remains to check the condition on $a_{ij}$:
\begin{align*}
     \sup\limits_{|x|\leq R} | a_{ij}^h(x) | 
     &\leq \sup\limits_{|x|\leq R}  h \,\Bigg|\,  \sum_{j=1}^N p_{ij}(x)\,(x_j - x_i)^2\, \Bigg| \\
     &\leq  h \, \sum_{j=1}^N 2R^2 \, \rightarrow 0 \text{ as } h\rightarrow 0.
\end{align*}
Hence both conditions of Lemma 8.5 are satisfied. This guarantees the weak convergence of the Markov process (which describes the ABM) to the ODE \eqref{Eqn: ODE model}.

\subsubsection{Variations}

We now consider changes to the way in which $i$ and $j$ are selected. If the selection is made by node degree \eqref{Eqn: node degree prob} then the only change is the replacement of the constant $\frac{1}{N^2}$ with $\frac{1}{N k_i}$ in $\Pi^h$. The application of Lemma 8.5 is then essentially the same as in the base case presented above, with this change in constants following through to give the normalisation in \eqref{Eqn: ODE model, k_i normalisation}. 

The second variation considers normalisation by likelihood of interaction, given by \eqref{Eqn: interaction prob normalisation}. However, individual $i$ then always interacts with the chosen individual $j$. This avoids the repetition of the appearance of $p_{ij}$. In this case $\Pi^h$ is given by
\begin{equation}
    \Pi^h(x,y) = 
    \begin{cases}
        \dfrac{1}{N}  p_{ij}(X) \bigg(\sum\limits_{l=1}^N p_{i l}(X) \bigg)^{-1} & \text{if } y = x +  e_i \mu^h (x_j - x_i) \text{ for some $i\neq j$}, \\[1em]
        \dfrac{1}{N} + \dfrac{1}{N}  \sum\limits_{i\neq j} \Bigg(1 -  p_{ij}(X) \bigg(\sum\limits_{l=1}^N p_{i l}(X) \bigg)^{-1} \Bigg) & \text{if } y = x,  \\[1.2em]
        0 & \text{otherwise.}  
    \end{cases}
\end{equation}
The change to the case $y=x$ does not affect directly the calculation of $a^h$ and $b^h$, so the only change that appears in the calculation is the replacement of $\frac{1}{N}$ with 
\begin{equation*}
    \frac{1}{N} \bigg(\sum\limits_{l=1}^N p_{i l}(X) \bigg)^{-1}.
\end{equation*}
As in the previous variation this change follows through to affect the normalisation only. The assumption that there exists $c>0$ such that $p_{ii}(X) > c$ for all $X\in\mathds{R}^N$ is required to ensure that the normalisation is always well-defined and the limiting coefficients remain Lipschitz continuous. 

It would be possible to consider the case that $i$ and $j$ go on to interact with probability $p_{ij}$. This would lead to the limiting ODE
\begin{equation} 
    \frac{dX_i}{dt} = \Bigg(\sum\limits_{l=1}^N p_{i l}(X)\Bigg)^{-1} \sum_{j=1}^N \big(p_{ij}(X)\big)^2\, (X_j - X_i).
\end{equation}
To the best of the authors' knowledge this system has not previously been studied as a model of opinion formation. However, it does appear as a limiting model when considering a fast adaptive network model in \cite{nugent2023evolving}.

We next consider a variation in which both individuals $i$ and $j$ update their opinions. Here the transition function is given by
\begin{equation}
    \Pi^h(x,y) = 
    \begin{cases}
        \dfrac{1}{N^2} \big(p_{ij}(x)+ p_{ji}(x)\big) & \text{if } y = x +  \mu^h \, (e_i - e_j) \,(x_j - x_i) \text{ for some $i< j$}, \\[1em]
        \dfrac{1}{N} + \dfrac{1}{N^2}  \sum\limits_{i\neq j} \big(1 - p_{ij}(x)\big) & \text{if } y = x,  \\[1.2em]
        0 & \text{otherwise.}  
    \end{cases}
\end{equation}
Note that we specify $i<j$ in the first case to avoid the repetition of states, but if we select a specific $i$ then interactions with any $j$ are possible due to the symmetry of the transition function. In this case,
\begin{align*}
    b^h_{i}(x) 
    &= \frac{1}{h} \sum_{j=1}^N \mu^h (x_j - x_i)\,\frac{1}{N^2} \big(p_{ij}(x)+p_{ji}(x)\big)\,, 
    = \frac{\mu^h}{N^2 h} \sum_{j=1}^N \big(p_{ij}(x)+p_{ji}(x)\big)\, (x_j - x_i)\,.
\end{align*}
Hence, in order to obtain again the same limiting $b_i$ we require that $p_{ij}(x)=p_{ji}(x)$ for all $x\in\mathds{R}^N$. Under this assumption we have
\begin{align*}
    b^h_{i}(x) &= \frac{2\mu^h}{N^2 h} \sum_{j=1}^N p_{ij}(x)\, (x_j - x_i)\,,
\end{align*}
 and so we require $\mu^h = Nh/2$. This also gives,
\begin{align*}
    a^h_{ii}(x) 
    &= \frac{h}{2} \Bigg( \sum_{j=1}^N p_{ij}(x)\,(x_j - x_i)^2\, \Bigg) \, \text{ and } \,
    a^h_{ij}(x) = -\frac{h}{2} p_{ij}(x)\,(x_j - x_i)^2 \,,
\end{align*}
so $a^h$ has the same limit at in the original ABM. Overall, under the assumption that $p_{ij}$ are symmetric, the only change required to obtain the same ODE limit \eqref{Eqn: ODE model} is multiplication by a constant in the value of $\mu^h$. 

\subsection{Additional sources of noise} \label{Appendix: Additional sources of noise}

Here we consider adaptations to the ABM described in Section \ref{Section: Additional sources of noise}. The general construction of the ABM is the same as in Section \ref{Subsection: ABM}, but with changes to the update rule \eqref{Eqn: Standard update scheme} to include additional noise terms.

The calculations of $a^h_{ij}$ and $b^h_{i}$ performed below are very similar to those above. As before, for an individual $i$, we sum over $j$ to include the effect of each possible interaction. When there is an additional source of noise we must also, for each $j$, integrate over the possible values the additional random variable could take. This is equivalent to integrating over the possible states that can be reached. 

\subsubsection{Ambiguity noise}

Consider the ABM with ambiguity noise, defined in Section \ref{Subsection: ambiguity noise} with update rule \eqref{Eqn: Ambiguity noise update scheme}:
\begin{equation*} 
    x_i(t+h) = 
    \begin{cases}
        x_i(t) + \mu^h\, \big(\omega_j(t) - x_i(t) \big) & \text{ with probability } \phi(|\omega_j - x_i|) \\
        x_i(t) & \text{ with probability } 1 - \phi(|\omega_j - x_i|) \,.
    \end{cases}
\end{equation*}
We will again use Lemma 8.5 from \cite{durrett2018stochastic} to show that this Markov process converges weakly to the ODE \eqref{Eqn: ODE model}.

Let $\eta$ be a family of random variables satisfying Assumption \ref{Assumption group: eta^h}. Let $\rho^h$ be the probability measure of the random variable $\eta^h$ and $\mathds{E}^h[\cdot]$ the expectation with respect to this measure. We obtain
\begin{align*}
    b^h_{i}(x) 
    &= \frac{1}{h} \sum_{j=1}^N \int_{\mathds{R}} \mu^h \, (x_j + \eta - x_i)\,\frac{1}{N^2} \phi\big(|x_j - x_i + \eta|\big) \,d\rho^h(\eta)  
    = \frac{1}{N} \sum_{j=1}^N \mathds{E}^h\big[(x_j - x_i + \eta)\, \phi\big(|x_j - x_i + \eta|\big) \big]   \\
    &= \frac{1}{N} \sum_{j=1}^N \Bigg( (x_j - x_i)\,\mathds{E}^h\big[ \phi\big(|x_j - x_i + \eta|\big) \big] + \mathds{E}\big[\eta\, \phi\big(|x_j - x_i + \eta|\big) \big] \Bigg)\,.  \\
    a^h_{ii}(x) 
    &= h \sum_{j=1}^N \int_{\mathds{R}} (x_j - x_i + \eta)^2\, \phi\big(|x_j - x_i + \eta|\big) \,d\rho^h(\eta)  \\
    &= h \sum_{j=1}^N \mathds{E}^h\big[(x_j - x_i + \eta)^2\, \phi\big(|x_j - x_i + \eta|\big) \big]   \\
    &= h \sum_{j=1}^N \Bigg( (x_j - x_i)^2\,\mathds{E}^h\big[ \phi\big(|x_j - x_i + \eta|\big) \big] + 2(x_j - x_i)\,\mathds{E}\big[\eta\, \phi\big(|x_j - x_i + \eta|\big) \big] + \mathds{E}^h\big[\eta^2\, \phi\big(|x_j - x_i + \eta|\big) \big] \Bigg)\,.  \\
\end{align*}
We now check that these approximations converge suitably to $b_i^h$ and $a_{ii} = 0$. Recall that $\phi$ is Lipschitz continuous with Lipschitz constant $L$ and is bounded above by $1$. Using this Lipschitz condition, together with triangle inequalities, we obtain,
\begin{align*}
    \sup\limits_{|x|\leq R}& | b_{i}^h(x) - b_{i}(x) | \\
    &\leq \sup\limits_{|x|\leq R} \frac{1}{N} \sum_{j=1}^N \bigg| (x_j - x_i)\,\mathds{E}^h\big[ \phi\big(|x_j - x_i + \eta|\big) \big] + \mathds{E}^h\big[\eta\, \phi\big(|x_j - x_i + \eta|\big) \big]  - (x_j - x_i)\, \phi\big(|x_j - x_i|\big) \bigg| \\
    &\leq \sup\limits_{|x|\leq R} \frac{1}{N} \sum_{j=1}^N 2R\, \Big| \mathds{E}^h\big[\phi\big(|x_j - x_i + \eta|\big) \big]  - \phi\big(|x_j - x_i|\big) \Big| 
    + \frac{1}{N} \sum_{j=1}^N \mathds{E}^h\big[\big| \eta\, \phi\big(|x_j - x_i + \eta|\big) \big| \big] \\
    &\leq \sup\limits_{|x|\leq R} \frac{1}{N} \sum_{j=1}^N 2R\, \mathds{E}^h\Big[ L \big| |x_j - x_i + \eta| - |x_j - x_i|\big| \Big]+ \mathds{E}^h\big[| \eta | \big] \\
    &\leq (2RL+1)\, \mathds{E}^h\big[ |\eta| \big] \,.
\end{align*}
As we assumed that $ \mathds{E}^h\big[ |\eta| \big] \rightarrow 0$ as $h\rightarrow0$, this gives the second condition of Lemma 8.5. The calculation for the bound on $a^h_{ii}$ is more straightforward as there is a factor of $h$, so we need only show that all other terms are bounded. 
\begin{align*}
    \sup\limits_{|x|\leq R} | a_{ii}^h(x) | 
    &\leq \sup\limits_{|x|\leq R} h \,\sum_{j=1}^N \Bigg( \Big| (x_j - x_i)^2\,\mathds{E}^h\big[ \phi\big(|x_j - x_i + \eta|\big) \big] \Big|+ \Big|2(x_j - x_i)\,\mathds{E}^h\big[\eta\, \phi\big(|x_j - x_i + \eta|\big) \big]\Big| \\
    &\quad \qquad \qquad + \Big|\mathds{E}^h\big[\eta^2\, \phi\big(|x_j - x_i + \eta|\big) \big] \Big|\Bigg) \\
    &\leq h \,\sum_{j=1}^N \Bigg( 2R^2 + 4R\, \mathds{E}^h[|\eta|] + \mathds{E}^h\big[\eta^2\big] \Big|\Bigg)\,.
\end{align*}
As $ \mathds{E}^h\big[ |\eta| \big] \rightarrow 0$ as $h\rightarrow0$, and $\mathds{E}^h\big[\eta^2\big]\leq C$ we have that the first condition of Lemma 8.5 is satisfied. 

\subsubsection{Adaptation noise}

Consider the ABM with adaptation noise, defined in Section \ref{Subsection: External noise} with update rule \eqref{Eqn: External noise with interactions update scheme}:
\begin{equation*} 
    x_i(t+h) = 
    \begin{cases}
        x_i(t) + \mu^h\, \big(x_j(t) - x_i(t) \big) + \xi^h & \text{ with probability } p_{ij}(x) \\
        x_i(t) & \text{ with probability } 1 - p_{ij}(x) \,.
    \end{cases}
\end{equation*}
We will now use Lemma 8.2 from \cite{durrett2018stochastic} to show that this Markov process converges weakly to the SDE \eqref{Eqn: limiting SDE external noise on interactions}. Note that this case is introduced after that of the ABM with external noise (using update rule \eqref{Eqn: External noise update scheme}) in Section \ref{Subsection: External noise}, but will be considered first here as the calculations for the external noise case are a simplification of those below. 

Let $\rho^h$ be the probability measure of the random variable $\xi^h$ and $\mathds{E}^h[\cdot]$ the expectation with respect to this measure. Assume the family of random variables $\xi$ satisfies Assumption \ref{Assumption group: xi^h}. We first consider the case in which noise is added only if individuals $i$ and $j$ interact. In this case,
\begin{align*}
    b^h_{i}(x) 
    &= \frac{1}{h} \sum_{j=1}^N \int_{\mathds{R}} \big(\mu^h \, (x_j - x_i) + \xi \big) \,\frac{1}{N^2} p_{ij}(x) \,d\rho^h(\xi)   \\
    &= \frac{1}{N} \sum_{j=1}^N p_{ij}(x) (x_j - x_i) \, + \frac{1}{N^2} \Bigg( \sum_{j=1}^N p_{ij}(x)\Bigg)\,\frac{\mathds{E}^h[\xi]}{h}  \,.  \\
    a^h_{ii}(x) 
    &= \frac{1}{h} \sum_{j=1}^N \int_{\mathds{R}} \big(\mu^h \, (x_j - x_i) + \xi \big)^2 \,\frac{1}{N^2} p_{ij}(x) \,d\rho^h(\xi)   \\
    &= h \sum_{j=1}^N p_{ij}(x) (x_j - x_i)^2\, + \frac{2}{N} \Bigg(\sum_{j=1}^N p_{ij}(x) (x_j - x_i)\Bigg)\, \mathds{E}^h[\xi] \, + \frac{1}{N^2} \Bigg(\sum_{j=1}^N p_{ij}(x) \Bigg)\, \frac{\mathds{E}^h[\xi^2]}{h} \,.
\end{align*}
As we wish to show that the limiting model is an SDE, rather than an ODE, we use Lemma 8.2 rather than Lemma 8.5. As before, we check the convergence of these approximations. 
\begin{align*}
    \sup\limits_{|x|\leq R} | b_{i}^h(x) - b_{i}(x) | &= \sup\limits_{|x|\leq R} \Bigg| \frac{1}{N^2} \Bigg( \sum_{j=1}^N p_{ij}(x)\Bigg)\,\frac{\mathds{E}^h[\xi]}{h} \Bigg| \, \leq \sup\limits_{|x|\leq R} \frac{1}{N} \Bigg| \frac{\mathds{E}^h[\xi]}{h} \Bigg|
\end{align*}
As we have assumed $m_1(\xi) = 0$, the first condition of Lemma 8.2 is satisfied. 
\begin{align*}
    \sup\limits_{|x|\leq R} | a_{ii}^h(x) - a_{ii}(x) | &=  
 \sup\limits_{|x|\leq R} \Bigg| h \sum_{j=1}^N p_{ij}(x) (x_j - x_i)^2\, + \frac{2}{N} \Bigg(\sum_{j=1}^N p_{ij}(x) (x_j - x_i)\Bigg)\, \mathds{E}^h[\xi] \, \\
    & \hspace{0.2\linewidth}+ \frac{1}{N^2} \Bigg(\sum_{j=1}^N p_{ij}(x) \Bigg)\, \frac{\mathds{E}^h[\xi^2]}{h} - \frac{m_2(\xi)}{N^2} \Bigg(\sum_{j=1}^N p_{ij}(x) \Bigg)\, \Bigg| \\
    &\leq h\big(4R^2\big)\, + \frac{4R}{N} \, \mathds{E}^h[\xi] \, + \frac{1}{N^2} \Bigg| \frac{\mathds{E}^h[\xi^2]}{h} - m_2(\xi) \Bigg|.
\end{align*}
Clearly the first term in this sum converges to $0$ as $h\rightarrow0$. As $m_1(\xi) = 0$, $\mathds{E}^h[\xi]$ also converges to $0$ as $h\rightarrow0$. In addition, the definition of $m_2(\xi)$ ensures that the last term in the sum also converges to $0$ as $h\rightarrow0$. Hence the second condition of Lemma 8.2 is satisfied.

As we are using Lemma 8.2, we also calculate $\gamma^h_p$ with $p=4$:
\begin{align*}
    \sup\limits_{|x|\leq R} \gamma^h_4(x) 
    &= \sup\limits_{|x|\leq R} \frac{1}{h} \sum_{i=1}^N \sum_{j=1}^N \int_{\mathds{R}} \big|\mu^h \, (x_j - x_i) + \xi \big|^4 \,\frac{1}{N^2} p_{ij}(x) \,d\rho^h(\xi) \,,  \\
    &\leq \frac{1}{N^2} \sum_{i=1}^N \sum_{j=1}^N p_{ij}(x) \Bigg( 16 N^4 h^3 R^4 + 32 N^3 h^2 R^3 \,\mathds{E}^h[\xi] + 24 N^2 h R^2 \,\mathds{E}^h[\xi^2] +8 N R \,\mathds{E}^h[\xi^3] + \frac{1}{h}\mathds{E}^h[\xi^4] \Bigg) \\
    &\leq h^3 \big(16 N^4 R^4 \big) + h^2 \big(32 N^3 R^3\big) \,\mathds{E}^h[\xi] + h \big(24 N^2 R^2 \big) \,\mathds{E}^h[\xi^2] +8 N R \,\mathds{E}^h[\xi^3] + \frac{1}{h}\mathds{E}^h[\xi^4] \,.
\end{align*}
The assumption that $m_3(\xi) = m_4(\xi) = 0$ guarantees that all terms in this sum converge to $0$ as $h\rightarrow0$. 

\subsubsection{External noise}

Now consider the ABM with external noise, defined in Section \ref{Subsection: External noise} with update rule \eqref{Eqn: External noise update scheme}:
\begin{equation*} 
    x_i(t+h) = 
    \begin{cases}
        x_i(t) + \mu^h\, \big(x_j(t) - x_i(t) \big) + \xi^h & \text{ with probability } p_{ij}(x) \\
        x_i(t) + \xi^h & \text{ with probability } 1 - p_{ij}(x) \,.
    \end{cases}
\end{equation*}
We will use Lemma 8.2 from \cite{durrett2018stochastic} to show that this Markov process converges weakly to the SDE \eqref{Eqn: limiting SDE external noise}. At each timestep, once individuals $i$ and $j$ are selected, there are two possible updates depending on whether or not an interaction occurs. This will lead to two integrals in the following calculation, but in fact gives a simpler final result. 

\begin{align*}
    b^h_{i}(x) 
    &= \frac{1}{h} \sum_{j=1}^N \int_{\mathds{R}} \big(\mu^h \, (x_j - x_i) + \xi \big) \,\frac{1}{N^2} p_{ij}(x) \,d\rho^h(\xi) \, + \frac{1}{h} \sum_{j=1}^N \int_{\mathds{R}} (\xi) \,\frac{1}{N^2} \big(1 - p_{ij}(x)\big) \,d\rho^h(\xi) \,,  \\
    &= \frac{1}{N} \sum_{j=1}^N p_{ij}(x) (x_j - x_i) \, + \frac{1}{N} \,\frac{\mathds{E}^h[\xi]}{h}  \,.  \\
    a^h_{ii}(x) 
    &= \frac{1}{h} \sum_{j=1}^N \int_{\mathds{R}} \big(\mu^h \, (x_j - x_i) + \xi \big)^2 \,\frac{1}{N^2} p_{ij}(x) \,d\rho^h(\xi) \, + \frac{1}{h} \sum_{j=1}^N \int_{\mathds{R}} (\xi)^2 \,\frac{1}{N^2} \big(1 - p_{ij}(x)\big) \,d\rho^h(\xi) \,,   \\
    &= h \sum_{j=1}^N p_{ij}(x) (x_j - x_i)^2\, + \frac{2}{N} \Bigg(\sum_{j=1}^N p_{ij}(x) (x_j - x_i)\Bigg)\, \mathds{E}^h[\xi] \, + \frac{1}{N} \frac{\mathds{E}^h[\xi^2]}{h} \,.
\end{align*}
Comparing these to $b^h_{i}(x)$ and $a^h_{ii}(x)$ from the previous case, the only change is the removal of the factor $\sum\limits_{j=1}^N p_{ij}(x)$ from the last term in each sum. Due to this similarity we do not show again the convergence of these approximations and $\gamma^h_p$, as the calculations are essentially the same as those above. 

\subsubsection{Random update distance}

Finally we consider the ABM with random update distances, defined in Section \ref{Subsection: Random update distance} with update rule \eqref{Eqn: Noisy update distance update scheme}:
\begin{equation*} 
    x_i(t+h) = 
    \begin{cases}
        x_i(t) + \nu^h\, \big(x_j(t) - x_i(t) \big) & \text{ with probability } p_{ij}(x) \\
        x_i(t) & \text{ with probability } 1 - p_{ij}(x) \,.
    \end{cases}
\end{equation*}
We will use Lemma 8.2 from \cite{durrett2018stochastic} to show that this Markov process converges weakly to the SDE \eqref{Eqn: limiting SDE noisy update distance}.

Let $\rho^h$ be the probability measure of the random variable $\nu^h$ and $\mathds{E}^h[\cdot]$ the expectation with respect to this measure. Assume the family of random variables $\nu$ satisfies Assumption \ref{Assumption group: nu^h}. Due to the linearity of the equations in $\mu$, the calculations in this section are much more straightforward: 
\begin{align*}
    b^h_{i}(x) 
    &= \frac{1}{h} \sum_{j=1}^N \int_{\mathds{R}} \nu \, (x_j - x_i) \,\frac{1}{N^2} p_{ij}(x) \,d\rho^h(\nu) 
    = \frac{\mathds{E}^h[\nu]}{h} \frac{1}{N^2} \sum_{j=1}^N p_{ij}(x) (x_j - x_i) \, .\\
    a^h_{ii}(x) 
    &= \frac{1}{h} \sum_{j=1}^N \int_{\mathds{R}} \nu^2 \, (x_j - x_i)^2 \,\frac{1}{N^2} p_{ij}(x) \,d\rho^h(\nu)  
    = \frac{\mathds{E}^h[\nu^2]}{h} \frac{1}{N^2} \sum_{j=1}^N p_{ij}(x) (x_j - x_i)^2 \,. 
\end{align*}
The conditions on $m_1(\nu)$ and $m_2(\nu)$ clearly ensure the first two conditions on Lemma 8.2 are satisfied. It only remains to calculate
\begin{align*}
    \sup\limits_{|x|\leq R} \gamma^h_4 
    &= \sup\limits_{|x|\leq R} \frac{1}{h} \sum_{i=1}^N \sum_{j=1}^N \int_{\mathds{R}} \nu^4 \, (x_j - x_i)^4 \,\frac{1}{N^2} p_{ij}(x) \,d\rho^h(\nu)  \\
    &= \frac{\mathds{E}^h[\nu^4]}{h} \frac{1}{N^2} \sum_{i=1}^N\sum_{j=1}^N p_{ij}(x) (x_j - x_i)^4   \\
    &\leq \frac{\mathds{E}^h[\nu^4]}{h} 16 R^4 \, .
\end{align*}
Hence the assumption that $m_4(\nu) = 0$ gives that the final condition of Lemma 8.2 is satisfied. 

\end{document}

%% file: preamble.tex
\usepackage[utf8]{inputenc}
\usepackage[dvipsnames]{xcolor}
\usepackage{amsfonts}
\usepackage{amsthm}
\usepackage{amsmath}
\usepackage{parskip}
\usepackage{mathrsfs}
\usepackage{amssymb}
\usepackage{xtab}
\usepackage{setspace}
\usepackage{datetime}
\usepackage{svg}
\usepackage{csquotes}
\usepackage{fancyhdr}
\usepackage{float}
\usepackage{comment}
\usepackage{enumitem}

\usepackage{amssymb,mathalfa,float,array,titlesec,mdframed,listings,amsmath,booktabs,dsfont,textcomp,graphicx,wrapfig,cases}
\usepackage[font=small,labelfont=bf]{caption}

\usepackage[utf8]{inputenc}
\usepackage[english]{babel}

\usepackage{float}

\usepackage{placeins}
\usepackage{subcaption}
\usepackage{hyperref}
\usepackage{cleveref}

\usepackage{cite}

\usepackage{multirow}

\graphicspath{{Figures/}}

\allowdisplaybreaks

\usepackage{thmtools}
\usepackage{thm-restate}

\declaretheorem[name=Definition,numberwithin=section]{definition}
\declaretheorem[name=Remark,numberwithin=section]{remark}

\declaretheorem[name=Assumption]{assumption}

\newcounter{subassumption}

%% file: main.bbl
\begin{thebibliography}{10}

\bibitem{hegselmann2004opinion}
Rainer Hegselmann.
\newblock Opinion dynamics: Insights by radically simplifying models.
\newblock {\em Laws and models in science}, pages 1--29, 2004.

\bibitem{deffuant2000mixing}
Guillaume Deffuant, David Neau, Frederic Amblard, and G{\'e}rard Weisbuch.
\newblock Mixing beliefs among interacting agents.
\newblock {\em Advances in Complex Systems}, 3(01n04):87--98, 2000.

\bibitem{lorenz2007continuous}
Jan Lorenz.
\newblock Continuous opinion dynamics under bounded confidence: A survey.
\newblock {\em International Journal of Modern Physics C}, 18(12):1819--1838,
  2007.

\bibitem{flache2017models}
Andreas Flache, Michael M{\"a}s, Thomas Feliciani, Edmund Chattoe-Brown,
  Guillaume Deffuant, Sylvie Huet, and Jan Lorenz.
\newblock Models of social influence: Towards the next frontiers.
\newblock {\em Journal of Artificial Societies and Social Simulation}, 20(4),
  2017.

\bibitem{sirbu2017opinion}
Alina S{\^\i}rbu, Vittorio Loreto, Vito~DP Servedio, and Francesca Tria.
\newblock Opinion dynamics: models, extensions and external effects.
\newblock {\em Participatory sensing, opinions and collective awareness}, pages
  363--401, 2017.

\bibitem{noorazar2020classical}
Hossein Noorazar, Kevin~R Vixie, Arghavan Talebanpour, and Yunfeng Hu.
\newblock From classical to modern opinion dynamics.
\newblock {\em International Journal of Modern Physics C}, 31(07):2050101,
  2020.

\bibitem{urbig2007communication}
Diemo Urbig and Jan Lorenz.
\newblock Communication regimes in opinion dynamics: Changing the number of
  communicating agents.
\newblock {\em arXiv preprint arXiv:0708.3334}, 2007.

\bibitem{laguna2004minorities}
Mar{\'\i}a~Fabiana Laguna, Guillermo Abramson, and Dami{\'a}n~H Zanette.
\newblock Minorities in a model for opinion formation.
\newblock {\em Complexity}, 9(4):31--36, 2004.

\bibitem{durrett2018stochastic}
Richard Durrett.
\newblock {\em Stochastic calculus: a practical introduction}.
\newblock CRC press, 2018.

\bibitem{fennell2021generalized}
Susan~C Fennell, Kevin Burke, Michael Quayle, and James~P Gleeson.
\newblock Generalized mean-field approximation for the deffuant opinion
  dynamics model on networks.
\newblock {\em Physical Review E}, 103(1):012314, 2021.

\bibitem{chu2022density}
Weiqi Chu and Mason~A Porter.
\newblock A density description of a bounded-confidence model of opinion
  dynamics on hypergraphs.
\newblock {\em arXiv preprint arXiv:2203.12189}, 2022.

\bibitem{ben2003bifurcations}
Eli Ben-Naim, Paul~L Krapivsky, and Sidney Redner.
\newblock Bifurcations and patterns in compromise processes.
\newblock {\em Physica D: nonlinear phenomena}, 183(3-4):190--204, 2003.

\bibitem{fortunato2005vector}
Santo Fortunato, Vito Latora, Alessandro Pluchino, and Andrea Rapisarda.
\newblock Vector opinion dynamics in a bounded confidence consensus model.
\newblock {\em International Journal of Modern Physics C}, 16(10):1535--1551,
  2005.

\bibitem{goddard2022noisy}
Benjamin~D Goddard, Beth Gooding, H~Short, and GA~Pavliotis.
\newblock Noisy bounded confidence models for opinion dynamics: the effect of
  boundary conditions on phase transitions.
\newblock {\em IMA Journal of Applied Mathematics}, 87(1):80--110, 2022.

\bibitem{como2009scaling}
Giacomo Como and Fabio Fagnani.
\newblock Scaling limits for continuous opinion dynamics systems.
\newblock In {\em 2009 47th Annual Allerton Conference on Communication,
  Control, and Computing (Allerton)}, pages 1562--1566. IEEE, 2009.

\bibitem{motsch2014heterophilious}
Sebastien Motsch and Eitan Tadmor.
\newblock Heterophilious dynamics enhances consensus.
\newblock {\em SIAM review}, 56(4):577--621, 2014.

\bibitem{ceragioli2021generalized}
Francesca Ceragioli, Paolo Frasca, Benedetto Piccoli, and Francesco Rossi.
\newblock Generalized solutions to opinion dynamics models with
  discontinuities.
\newblock In {\em Crowd Dynamics, Volume 3: Modeling and Social Applications in
  the Time of COVID-19}, pages 11--47. Springer, 2021.

\bibitem{koponen2022agent}
Ismo~T Koponen.
\newblock Agent-based modeling of consensus group formation with complex webs
  of beliefs.
\newblock {\em Systems}, 10(6):212, 2022.

\bibitem{steiglechner2023noise}
Peter Steiglechner, Marijn~A Keijzer, Paul~E Smaldino, Deyshawn Moser, and
  Agostino Merico.
\newblock Noise and opinion dynamics: How ambiguity promotes pro-majority
  consensus in the presence of confirmation bias.
\newblock {\em SocArXiv}, 2023.

\bibitem{ceragioli2012continuous}
Francesca Ceragioli and Paolo Frasca.
\newblock Continuous and discontinuous opinion dynamics with bounded
  confidence.
\newblock {\em Nonlinear Analysis: Real World Applications}, 13(3):1239--1251,
  2012.

\bibitem{brooks2022emergence}
Heather~Z Brooks, Philip~S Chodrow, and Mason~A Porter.
\newblock Emergence of polarization in a sigmoidal bounded-confidence model of
  opinion dynamics.
\newblock {\em arXiv preprint arXiv:2209.07004}, 2022.

\bibitem{blondel2010continuous}
Vincent~D Blondel, Julien~M Hendrickx, and John~N Tsitsiklis.
\newblock Continuous-time average-preserving opinion dynamics with
  opinion-dependent communications.
\newblock {\em SIAM Journal on Control and Optimization}, 48(8):5214--5240,
  2010.

\bibitem{lacker2018mean}
Daniel Lacker.
\newblock Mean field games and interacting particle systems.
\newblock {\em preprint}, 2018.

\bibitem{nugent2023evolving}
Andrew~J Nugent, Susana~N Gomes, and Marie-Therese Wolfram.
\newblock On evolving network models and their influence on opinion formation.
\newblock {\em arXiv preprint arXiv:2305.09483}, 2023.

\bibitem{boghosian2022particle}
Bruce Boghosian, Christoph B{\"o}rgers, Natasa Dragovic, Anna Haensch, and
  Arkadz Kirshtein.
\newblock A particle method for continuous hegselmann-krause opinion dynamics.
\newblock {\em arXiv preprint arXiv:2211.06265}, 2022.

\bibitem{erdHos1960evolution}
Paul Erd{\H{o}}s, Alfr{\'e}d R{\'e}nyi, et~al.
\newblock On the evolution of random graphs.
\newblock {\em Publ. Math. Inst. Hung. Acad. Sci}, 5(1):17--60, 1960.

\bibitem{pineda2009noisy}
Miguel Pineda, Raul Toral, and Emilio Hernandez-Garcia.
\newblock Noisy continuous-opinion dynamics.
\newblock {\em Journal of Statistical Mechanics: Theory and Experiment},
  2009(08):P08001, 2009.

\bibitem{billingsley2013convergence}
Patrick Billingsley.
\newblock {\em Convergence of probability measures}.
\newblock John Wiley \& Sons, 2013.

\end{thebibliography}
